\documentclass[reqno,11pt]{amsart}
 
 \RequirePackage{amsthm,amsmath,amsfonts,amssymb}
 
 \usepackage{mathtools}
 
\numberwithin{equation}{section}
 
  \theoremstyle{plain} 
 \newtheorem{theorem}{Theorem}[section]
\newtheorem{lemma}[theorem]{Lemma}
\newtheorem{corollary}[theorem]{Corollary}

\theoremstyle{remark} 
\theoremstyle{definition}
\newtheorem{assumption}[theorem]{Assumption}
\newtheorem{definition}[theorem]{Definition}
\newtheorem{remark}[theorem]{Remark}
 
 \makeatletter 
 \def\dashint{\operatorname{\,\,\,\mathclap{\int}\!\!\!\!\!\text{\phantom{n}\bf  --}\!\!}}
 \makeatother
 
\makeatletter 
\def\dashnorm{\,\,\text{\bf--}\kern-.5em\|}
\def\ninf{\qopname\relax\@empty{inf\phantom{p}\!\!\!}}

 \makeatother

\newcommand\sfp{{\sf p}}
\newcommand\sfq{{\sf q}}

\newcommand\bB{\mathbb{B}}
\newcommand\bC{\mathbb{C}}

\newcommand\bR{\mathbb{R}}
\newcommand\bS{\mathbb{S}}

\newcommand \bW{\mathbb{W}}
\newcommand\bZ{\mathbb{Z}}

\newcommand\cF{\mathcal{F}}

\newcommand \cL{\mathcal{L}}

\newcommand\cN{\mathcal{N}}

\newcommand\frp{\mathfrak{p}}
\newcommand\frq{\mathfrak{q}}

\newcommand{\loc}{{\rm loc}\,}

 \newcommand{\mysection}[1]{\section{#1}
 \setcounter{equation}{0}}
\newcommand\cSigma{{\scriptstyle \Sigma}}

\newcommand{\nliminf}{\operatornamewithlimits{\underline{lim}}}

\begin{document}

 \title[On strong solutions]
{On strong solutions of time inhomogeneous
It\^o's equations with Morrey diffusion 
gradient and drift. A supercritical case}
\author{N.V. Krylov}

\email{nkrylov@umn.edu}
\address{School of Mathematics, University of Minnesota, Minneapolis, MN, 55455}
 
 

\begin{abstract} 
We prove the existence of strong solutions
of It\^o's stochastic time dependent equations with
irregular diffusion and drift terms
of Morrey  spaces. Strong uniqueness is also
discussed.
\end{abstract}

\maketitle

\mysection{Introduction}
                                                  \label{section 3.11.1}

This paper is a natural continuation of
 \cite{1}, where the strong existence and strong uniqueness is proved for the time {\em homogeneous\/} stochastic equations with Morrey drift. Here we do the same
in the time inhomogeneous case and obtain the results
containing those in \cite{1}.

Let $\bR^{d}$ be a $d-$dimensional Euclidean space of points
$x=(x^{1},...,x^{d})$ with $d\geq 3$. Let $(\Omega,\cF,P)$ be a 
complete probability space,
carrying a $d_{1} $-dimensional Wiener process
  $w_{t}$ with $d_{1}\geq d$. Fix $\delta\in(0,1]$ and denote by $\bS_{\delta}$ the set of $d\times d$ symmetric matrices
whose eigenvalues lie in $[\delta,\delta^{-1}]$.

Assume that on $\bR^{d+1}=\{(t,x):t\in\bR,x\in\bR^{d}\}$ we are given Borel $\bR^{d}$-valued 
functions  $b=(b^{i})$ and $\sigma^{k}=(\sigma^{ik})$,
$k=1,...,d_{1}$. Set $\sigma=(\sigma^{ik})$.
We are going
to   investigate
the equation
\begin{equation}
                                         \label{6.15.2}
 x_{s}=x +\int_{0}^{s}\sigma  (t+u,x_{u})\,dw  _{u}+
\int_{0}^{s}b(t+u,x_{u})\,du.
\end{equation}
The basic assumption is that with $a^{ij}=\sigma^{ik}
\sigma^{jk}$ we have for all values of the arguments that
$$
a:=(a^{ij})\in \bS_{\delta}.
$$

We are interested in the so-called strong solutions, that is
solutions such that, for each $s\geq0$, $x_{s}$ is $\cF^{w}_{s}$-measurable,
where $\cF^{w}_{s}$ is the completion of $\sigma(w_{u}:u\leq s)$. We present
sufficient conditions for the equation to have a strong solution  and deal with its uniqueness.

After the classical work by K. It\^o showing that there exists
a unique strong solution of \eqref{6.15.2} if $\sigma $ and $b$
are Lipschitz continuous in $x$ (may also depend on   $\omega$),
much effort was applied to relax these Lipschitz conditions.
In   case $d=d_{1}=1$ T. Yamada and S. Watanabe \cite{YW_71} relaxed
the Lipschitz condition on $\sigma$ to the H\"older $(1/2)$-condition
(and even slightly weaker condition) and kept $b$ Lipschitz
(slightly less restrictive). Much attention was paid to equations
with continuous coefficients
satisfying the so-called monotonicity conditions
(see, for instance, \cite{Kr_84} and the references therein).

T. Yamada and S. Watanabe \cite{YW_71} also put forward
a very strong theorem, basically, saying that the existence of weak solutions and strong uniqueness
implies the existence of strong
solutions. Unlike the present paper, 
the majority of papers on the subject after that time
are using their theorem.
 S. Nakao (\cite{Na_72}) proved
the strong solvability in time homogeneous case
 if $d=d_{1}=1$ and $\sigma$ is bounded away from zero and infinity
 and is
locally of bounded variation. He also assumed that $b$ is bounded,
but from his arguments it is clear that the summability of $|b|$
suffices. In this respect his result basically shows that
our results are also true if $d=1$
and the coefficients are independent of time. However, the 
general case that $d=2$ is quite open.
  By the way, the restriction $d\geq3$
appears because 
the method of proof
of Lemma \ref{lemma 4.28.1} requires $2\leq  p_{0}<  p_{b}\leq d$
for certain $  p_{0}, p_{b}$.

A. Veretennikov was the first author who in \cite{Ve_80} 
not only proved the existence of strong solutions
in the time inhomogeneous {\em multidimensional\/} case  when $b$ is bounded,
but also considered the case of $\sigma $ in Sobolev class,
namely,   $\sigma _{x}\in L_{2d,\loc}$. He used A. Zvonkin's method (see
\cite{Zv_74}) of transforming the equation in such a way that the
drift term disappears.
In \cite{XXZZ_20} (also see the references
there) the result of Veretennikov is
extended to the case of $\sigma$ uniformly continuous in $x$ and $\sigma_{x},b\in L_{p,q}$
with, perhaps, different $p,q$ for $\sigma_{x}$ and $b$ satisfying
\begin{equation}
                           \label{6.12.1}
\frac{d}{p}+\frac{2}{q}< 1
\end{equation}
(so-called subcritical Ladyzhenskaya-Prodi-Serrin
condition).
In that case much information is available,
we refer the reader
to \cite{Zh_11},  \cite{Zh_20}, \cite{XXZZ_20}, and the references therein.  

Even the case when   $\sigma $  is constant and the process is nondegenerate
attracted
very much attention.   
M. R\"ockner and the author in \cite{KR_05} proved,
among other things, the existence of strong solutions 
when $b\in L_{p,q}$ under condition \eqref{6.12.1}.
We  refer to \cite{PFPR_13}, \cite{MNF_15},  \cite{Zh_16} and the references therein for further results in  this direction.
If $b$ is bounded  A.~Shaposhnikov (\cite{Sh_14}, \cite{Sh_17})
proved the so called path-by-path uniqueness, which, basically,
means that for almost any trajectory $w_{t}$ there is only one
solution (adapted or not). This result was already announced
by A. Davie before with a very entangled proof which left many doubtful.

In a fundamental work by L.~Beck, F.~Flandoli, M.~Gubinelli, and M.~Maurelli
(\cite{BFGM_19}) the authors investigate such equations from
 the points of view of It\^o stochastic equations, stochastic transport 
equations, and stochastic continuity equations. Their article
contains an enormous amount of information and a vast references list. 
In what concerns our situation they require ($\sigma $ constant
and the process is nondegenerate) what they call LPS-condition (slightly imprecise): $b\in L_{p,q,\loc}$,
$q<\infty$, with equality in \eqref{6.12.1} in place of $<$,
or $p=d$ but $\|b\|_{L_{p,\infty}}$ to be sufficiently small, or else that $b(t,\cdot)$ to
be continuous as an $L_{d}(\bR^{d})$-function,
and they prove strong solvability and strong uniqueness
(actually, path-by-path-uniqueness which is stronger) but only for {\em almost
all\/} starting points $x$. 
Concerning the strong solutions
starting from any point $x$
in the {\em time dependent\/} case
with singular $b$  and constant $\sigma$ probably the best results to date belong to R\"ockner and Zhao \cite{RZ_21}, where, among very many other things, they prove existence and uniqueness of strong solutions of
equations like \eqref{6.15.2} with $b\in L_{p,q}$ and $p,q<\infty$, 
with equality in \eqref{6.12.1} in place of $<$, or when $b(t,\cdot)$
is continuous as an $L_{d}(\bR^{d})$-function.
  
We refer the reader to
\cite{BFGM_19} and \cite{RZ_21} also for
a very good review  of the motivation
related to the Navier-Stokes equation and history of the problem.

Our approach is absolutely different from
all articles mentioned above and all articles
which one can find in their references.
 We do not use
Yamada-Watanabe theorem or transformations of the noise or   a compactness criterion for random
fields in Wiener-Sobolev spaces as in \cite{RZ_21}.
Instead, our method is inspired by
an analytic criterion for the existence
of strong solutions which first appeared in \cite{VK_76},
 some 45 years ago and was first used only
in \cite{Kr_21}.
   To make this method work we use  ideas
from many papers, most relevant of which are  the results in
\cite{DK_18}, \cite{Kr_21}, \cite{1}, \cite{10}, \cite{13}, \cite{VK_76}.
 
Our results are not covered by 
\cite{RZ_21} not only because we have
variable $\sigma  $  but also because,
for any $t$, our
$b(t,\cdot)$ may be generally not even in $L_{d}$.

A very rich literature
on the {\em weak\/} solutions of \eqref{6.15.2}     is beyond
  the scope of this article.  By the way,
G. Zhao (\cite{Zh_20_1}) gave an example showing
that, if in condition \eqref{3.14.2}   we replace $r$ with $r^{ \alpha}$, $\alpha>1$, the weak uniqueness
may fail even in the time homogeneous case
and unit diffusion.

Here is an example in which we prove existence
and conditional uniqueness of strong solutions. Take $d=3,d_{1}=12$, and
for some numbers $\alpha,\beta,\gamma\geq0$
let $\sigma^{k}$ be the $k$th column in
 ($0/0:=3^{-1/2}$)
\begin{equation}
                                    \label{6.3.4}
\begin{pmatrix}\alpha & 0 & 0\\
0 & \alpha & 0\\
0 & 0 & \alpha
\end{pmatrix},
\frac{\beta}{|x|}
\begin{pmatrix}
x^{1} &  x^{2} & x^{3} & 0 &  0   & 0   & 0  & 0  & 0 \\
0 & 0 & 0 & x^{1}  & x^{2}   & x^{3}   & 0  & 0  & 0 \\
0 & 0 & 0 & 0  & 0   & 0  & x^{1}& x^{2} & x^{3} 
\end{pmatrix} ,
\end{equation}
$$
b(x)=-\frac{\gamma}{|x|}\,\frac{x}{|x|}I_{0<|x|\leq 1} +\xi(t)\eta(t,x),
$$
 where $\eta$ is bounded 
$\bR^{3}$-valued  and $\xi$ is real-valued of class $
 L_{2}(\bR)$.
Our result shows that if $\alpha=1$ and
$\beta$ and $\gamma$ are sufficiently small,
then \eqref{6.15.2} has a strong solution
which is conditionally unique, however, if
$\xi\equiv0$, then any solution is strong and unique.
By the way, if $\xi\equiv0$, $\alpha=\gamma=0$ and $\beta=1$, there exist
strong solutions of \eqref{6.15.2}
only if the starting point $x\ne0$ (see \cite{Kr_21}). In case $\alpha=1$ and $\beta=0$
strong solutions exist only if $\gamma$ is sufficiently small.  In case   $\alpha=1,\beta=\gamma=0$,
absent in \cite{RZ_21}, the authors of 
\cite{BFGM_19} prove the unique strong solvability
only for almost all starting points.

Observe that for $\beta\ne 0$
and $\gamma\ne0$ we have $D\sigma ,b(t,\cdot)\in L_{d-\varepsilon,\loc}(\bR^{d})$ for any $\varepsilon\in(0,1)$
but not for $\varepsilon=0$. Recall  that the case
of time independent $\sigma,b$ with $D\sigma ,b\in L_{d,\loc}$ is investigated
in \cite{Kr_21} the main idea of which
is used here as well.

Other examples can be found in Remark 
\ref{remark 6.4.4}. There as above we compare
our results with the ones obtained when $\sigma$
is the unit matrix. In this connection note that
our results are new even if $b\equiv0$.

We conclude the introduction by some notation.
We set  
$$
 D_{i}u =\frac{\partial u}{\partial x^{i}} ,\quad Du =(D_{i}u),\quad D_{ij}u =D_{i}D_{j}u,\quad
D^{2}u=(D_{ij}u ),\quad
\partial_{t}u=\frac{\partial u}{\partial t}.
$$

 If $\sigma =(\sigma^{ij...})$ by $|\sigma|^{2}$
we mean the sum of squares of all entries.

For $p\in[1,\infty)$, 
and domain $\Gamma\subset \bR^{d}$ by   $L_{p}(\Gamma)$  we mean the space 
of Borel (real-, vector-, matrix-valued...)
 functions on   $\Gamma $  with finite norm given by
$$
 \|f\|_{L_{p}(\Gamma)}^{p}=\int_{\Gamma}|f(x)|^{p}\,dx .
$$
Set $L_{p}=L_{p}(\bR^{d})$.

For $p,q\in[1,\infty)$ and domain $Q\subset\bR^{d+1}$ by $L_{p,q}(Q)$
we mean the space of Borel (real-, vector- or matrix-valued)
 functions on $Q$   with finite norm given by
\begin{equation}
                                       \label{6.4.1}
\|f\|_{L_{p,q}(Q)}^{q}=\|fI_{Q}\|_{L_{p,q}}^{q}
=\int_{\bR}\Big(\int_{\bR^{d}}|fI_{Q}(t,x)|^{p}\,
dx\Big)^{q/p}\,dt.
\end{equation}
Set $L_{p,q}=L_{p,q}(\bR^{d+1})$,
$$
\dashnorm f\|_{L_{p,q}(Q)}=\|1\|^{-1}_{L_{p,q}(Q)}\|f\|_{L_{p,q}(Q)}.
$$

By $W^{2}_{p}$ we mean the space 
of Borel functions $u$ on $\bR^{d}$ whose Sobolev derivatives 
$Du$ and $D^{2}u$ exist and $u,Du ,D^{2}u \in L_{p}$.
The  norm in $W^{2}_{p}$ is given by
$$
\|u\|_{W^{2}_{p}}=\|D^{2}u \|_{L_{p}}+\|u \|_{L_{p}}.
$$ 
 As usual, we write $f\in L_{p,\loc}$
if $f\zeta\in L_{p}$ for any $\zeta\in C^{\infty}_{0}$ ($=C^{\infty}_{0}
(\bR^{d})$).  

By $W^{1,2}_{p,q}(Q)$ we mean the collection
of $u$ such that $\partial_{t}u$, $D^{2}u $, $Du $, $u
\in L_{p,q}(Q)$. The norm in $W^{1,2}_{p}(Q)$
is introduced in an obvious way.
We abbreviate $W^{1,2}_{p,q } =W^{1,2}_{p,q}(\bR^{d+1})$.

If a Borel $\Gamma\subset \bR^{d+1}$, by $|\Gamma|$ we mean its Lebesgue
measure,
$$
\dashint_{\Gamma}f(t,x)\,dxdt=\frac{1}{|\Gamma|}
\int_{\Gamma}f(t,x)\,dxdt.
$$
Similar notation is used for $\Gamma\subset
\bR^{d}$.

  Introduce
$$
B_{R}(x)=\{y\in\bR^{d}:|x-y|<R\},\quad B_{R} =B_{R}(0),
$$
and let $\bB_{R}$ be the collection of balls of radius $R$. Also let  
$$
C_{\tau,\rho}(t,x)=[t,t+\tau)\times B_{\rho}(x),\quad C_{\rho}...=C_{\rho^{2},\rho}...,\quad C_{\rho}=C_{\rho}(0,0),
$$
and let $\bC_{\rho}$ be the collection of
$C_{\rho}(t,x)$.

For    $\beta\geq 0$, $p,q\in[1,\infty)$ 
and domain $Q\subset \bR^{d+1}$ introduce
Morrey's space $E_{p,q,\beta}(Q) $
as the set of $g\in  L_{p,q,\loc}(Q)$ such that  
\begin{equation}
                             \label{8.11.02}
\|g\|_{E_{p,q,\beta}(Q) }:=
\sup_{\rho\leq 1,C\in\bC_{\rho}}\rho^{\beta}
\dashnorm gI_{Q}  \|_{ L_{p,q}(C)} <\infty .
\end{equation}

Define
$$
E^{1,2}_{p,q,\beta}(Q)  =\{u:u,Du,D^{2}u,
\partial_{t}u\in E_{p,q,\beta}(Q)  \}
$$
and provide $E^{1,2}_{p,q,\beta}(Q)  $ with an obvious norm. If $Q=\bR^{d+1}$, we drop $Q$ in the above notation.

It is important to have in mind that if
$\beta<2$ (our main case) and $u\in E^{1,2}_{p,q,\beta}$,
then, according to Lemma 2.5 of \cite{11},
$u$ is bounded and continuous. 

In the proofs of our results
we use various (finite) constants called $N$ which
may change from one occurrence to another
and depend on the data only in the same way
  indicated in the statements
of the results unless specifically pointed out
otherwise.

\mysection{Setting and main results}
 
                 \label{section 12.20.1}

Fix some numbers $p_{D\sigma} \in(2,d], p_{b}\in (2\vee(d/2), d] $ (hence $d\geq3$)
and $\rho_{D\sigma},\rho_{b}\in(0,\infty)$.

\begin{definition}
                  \label{definition 12.16.1}
We call a real- or vector- or else tensor-valued
function $f(t,x)$ given on $ 
 \bR^{d+1} $ {\em admissible\/} if it  is
represented as $f = f_{M} +f_{B} $ (``Morrey part'' of
$f$ plus its ``bounded part'')
with  Borel $f_{M}$ and $f_{B}$   such that there exists   numbers $p_{f}\geq1$, $\rho_{f}\in(0,\infty)$, and  a {\em constant\/} $\hat f_{M}<\infty$ for which
\begin{equation}
                           \label{7.24.4}
\Big(\dashint_{B } |f_{M}(t,x)|^{  p_{f}}\,dx\Big)^{1/ p_{f}}\leq \hat f_{M}  \rho^{-1},
\end{equation}
whenever  $t in\bR $,  $B\in\bB_{\rho}$ and $\rho\leq \rho_{f}$,
and there exists a measurable $\bar f (t)\geq
\sup_{x\in\bR^{d}}|f_{B} ( t,x)|$
such that  
\begin{equation}
                           \label{7.26.1}
\int_{\bR}|\bar f (t) | ^{2}\,dt <\infty.
\end{equation}
 
\end{definition}

To be more concrete by $\hat f_{M}$ we mean the least
constant for which \eqref{7.24.4} holds
whenever  $t \in\bR $,  $B\in\bB_{\rho}$ and $\rho\leq \rho_{f}$.
Introduce
$$
\beta_{f}(t)=\sup_{s}\int_{s}^{s+t}
\bar f^{2}(r)    \,dr, \quad t\geq0.
$$
 
Assume that the Sobolev derivatives $D\sigma$ exist
and recall that $a=\sigma\sigma^{*}$.
\begin{assumption}
                  \label{assumption 12.19.1}
$  D\sigma$ and  $b  $ are admissible.

\end{assumption}

\begin{remark}
                            \label{remark 7.21.2}
That $D\sigma$ is admissible implies that $Da$ is also admissible with
$Da_{M}:=(D\sigma_{M})\sigma^{*}+\sigma D\sigma^{*}_{M}$, $p_{Da}=p_{D\sigma}$, $\widehat{Da_{M}}
\leq N(d,d_{1},\delta)\widehat{D\sigma_{M}}
$, and $\beta_{Da}\leq N(d,d_{1},\delta)\beta_{D\sigma}$.
\end{remark}

Fix some $p_{0},q_{0}$ satisfying
\begin{equation}
                         \label{5.31.2}
p_{0}\in (2 , p_{b}),\quad q_{0}\in(1,\infty),
\quad  \frac{d}{p_{0}}+\frac{2}{q_{0}}<2.
\end{equation}
 Observe that the last inequality in \eqref{5.31.2}
implies that $p_{0},p_{b}>d/2$. 

\begin{theorem}
                        \label{theorem 7.18.1}
Under Assumption \ref{assumption 12.19.1}
there exists $\varepsilon_{0},\varepsilon_{0}'>0$ depending only
on $d,d_{1},\delta,   p_{b},p_{D\sigma},p_{0},q_{0}$,  with
$\varepsilon_{0}'$ also depending on
  $\rho_{D\sigma}$, such that, if
  $\widehat{{D\sigma_{M}}}\leq \varepsilon_{0},
 \hat b_{M}\leq \varepsilon'_{0}$, then for
each $(t,x)$ 
equation \eqref{6.15.2} has a strong
solution  possessing the property 

(a)  for
 any $T\in(0,\infty)$, $m=1,2,...$,
there exists a constant $N$
such that for any Borel nonnegative $f$
\begin{equation}
                                \label{5.31.1}
E \Big(\int_{0}^{T}f( s,x_{s})\,ds\Big)^{m}\leq 
N \|f\|^{m}_{L_{p_{0},q_{0} }}.
\end{equation}
 
Furthermore (conditional uniqueness), if there exists a  solution $y_{s}$ of \eqref{6.15.2} with the same initial condition having the property 

($\text{a}\,'$) for each $T\in(0,\infty)$
there is a constant $N$ such that for any Borel nonnegative $f$
\begin{equation}
                                \label{7.19.1}
E  \int_{0}^{T}f( s,y_{s})\,ds \leq 
N \|f\| _{L_{p_{0},q_{0}}},
\end{equation}
then $x_{\cdot}=y_{\cdot}$ (a.s.).
\end{theorem}

This theorem is proved in Section \ref{section 7.29.1}.
It is convenient to prove it
by using a somewhat weaker result (Theorem \ref{theorem 5.30.1} below) 
 imposing an additional assumption on $b$.

\begin{assumption}
                      \label{assumption 3.14.1} 
  For some $r_{b}\in(0,1]$,  $\frp_{b},\frq_{b}\in (1,\infty)$,
 such that
\begin{equation}
                           \label{3.14.1}
 \frac{d}{\frp_{b}}
+\frac{2}{\frq_{b}} \geq 1 ,
 \end{equation}
we have 
\begin{equation}
                           \label{3.14.2}
 \sup_{r\leq r_{b}}r
\sup_{C\in \bC_{r}} 
\dashnorm b \|_{L_{\frp_{b},\frq_{b}}(C)}= \hat b <\infty.
\end{equation}

\end{assumption}

Assumptions \ref{assumption 12.19.1} and \ref{assumption 3.14.1} are supposed to be satisfied
below and
throughout Sections \ref{section 6.4.1} and \ref{section 6.4.2} with the restrictions
on $\widehat{D\sigma_{M}},\hat b_{M},\hat b$ and
$\beta_{D\sigma}(\rho_{\sigma}^{2})$ meaning that
they are ``small enough''  
made precise in somewhat technical Remark \ref{remark 7.21.1}
and Assumption \ref{assumption 7.21.2},  
which is also supposed to be satisfied
below and
throughout Sections \ref{section 6.4.1} and \ref{section 6.4.2}. 

\begin{theorem}[unconditional strong uniqueness]
                         \label{theorem 7.24.1}
Under Assumptions \ref{assumption 12.19.1},
\ref{assumption 3.14.1},
and 
\ref{assumption 7.21.2}  suppose that
$$
\frac{d} {\frp_{b}}+\frac{1}{\frq_{b}}\leq 1.
$$
Also suppose that, if  $\frp_{b}\geq \frq_{b}$,   then $\hat b$ is defined
as in \eqref{3.14.2}, but if $\frp_{b}\leq \frq_{b}$,
$\dashnorm b \|_{L_{\frp_{b},\frq_{b}}(C)}$  in \eqref{3.14.2} is
understood as
\begin{equation}
                                     \label{7.29.6}
\|I_{C}\|_{L_{\frp_{b},\frq_{b}}}^{-1}\Big(\int_{\bR^{d}}
\Big(\int_{\bR}I_{C}|b| ^{\frq_{b}}\,dt\Big)^{\frp_{b}/\frq_{b}}
\,dx\Big)^{1/\frp_{b}}.
\end{equation}
Then there exists one and only one strong solution of \eqref{6.15.2} for each initial data $(t,x)$
and any solution is strong.

\end{theorem}

Proof. By Theorem \ref{theorem 7.18.1} there exists
a strong solution possessing property (a).
By Theorem 3.8 of \cite{13} weak uniqueness holds.
After that it only remains to use Theorem 5.1 of \cite{10}. The theorem is proved. \qed

Just in case, note that, in the proof of Theorem
\ref{theorem 7.18.1}, Theorem \ref{theorem 7.24.1}
is not used.

\begin{theorem}
                             \label{theorem 5.30.1}
     Let $p_{0}<\frp_{b},q_{0}<\frq_{b}$. Then 
 under Assumptions \ref{assumption 12.19.1},
\ref{assumption 3.14.1}, and 
\ref{assumption 7.21.2} 
 equation \eqref{6.15.2} has a strong
solution $x_{\cdot}$ (for    $(t,x)=(0,0)$ )  possessing  property
(a). Furthermore, any solution starting from the origin and possessing property ($\text{a}\,'$) is strong
and coincides with~$x_{\cdot}$.

\end{theorem}

The proof of this theorem is given in Section \ref{section 6.4.2}.

\begin{theorem}
                            \label{theorem 6.1.1}
Under the assumptions of Theorem 
\ref{theorem 5.30.1} the strong solution starting from $(t,x)$
is also such that, for
any $T\in(0,\infty)$ and $p,q$ satisfying
\begin{equation}
                         \label{6.1.2}
 p\in[p_{0},\infty),\quad q\in(1,\infty),
\quad p_{0}\Big(\frac{1}{\alpha}-\frac{1}{q_{0}}\Big)
< p\Big(1-\frac{1}{q }\Big),
\end{equation}
where $\alpha=\alpha(p_{0},q_{0})>1$ is specified in the proof,
there exists a constant $N$
such that for any Borel nonnegative $f$
\begin{equation} 
                                \label{6.1.1}
E \int_{0}^{T}f( s,x_{s})\,ds\leq 
N\Big(\int_{0}^{T}\Big(\int_{\bR^{d}} \|f(s,\cdot)\|_{L_{p }(B_{1}(z))}^{p}\Psi(s,z-x)\,dz\Big)^{q/p}ds\Big)^{1/q},
\end{equation}
where $\Psi(s,y)=I_{B_{1}}+e^{-|y| /(Ns)}$.
\end{theorem}

The proof of this theorem is given in Section \ref{section 6.4.2}.

\begin{remark}
                             \label{remark 7.29.1}

Condition \eqref{6.1.2} allows to $q$ be close to
1 on the account of taking $p$ large enough
and even for $p=p_{0}$, $q=q_{0}$ estimate
\eqref{6.1.1} is stronger than \eqref{7.19.1}.
In this case, as it follows from the proof,
one can take $q$ slightly smaller than $q_{0}$.
\end{remark}

\begin{remark}
                      \label{remark 6.2.1}
Probably the best prior result about existence
of strong solutions for It\^o equations
with singular drift belong to R\"ockner and Zhao
\cite{RZ_21}. In this paper they consider the case
when $\sigma$ is {\em the unit matrix\/} and 
there exist $p,q$ such that
\begin{equation}
                                   \label{6.2.2}
b\in L_{p,q},\quad
p,q\in (1,\infty),\quad \frac{d}{p}+\frac{2}{q}=1
\end{equation}
or $b\in C([0,T],L_{d})$, and they prove that, for any initial
data, equation \eqref{6.15.2} has a strong solution
on $[0,T]$
possessing the property

(b) given any $p,q$ satisfying
\begin{equation}
                                \label{7.29.5}
p,q\in(1,\infty),\quad \frac{d}{p}+\frac{2}{q}<2
\end{equation}
it holds that
\begin{equation}
                                \label{7.29.4}
E  \int_{0}^{T}|f( s,x_{s})|\,ds \leq 
N \sup_{C\in \bC_{1}}\|f \|_{L_{p,q}(C)},
\end{equation}
where $N$ is independent of $f$.

They also prove
that  {\em strong\/} solutions, possessing 
  property (b),  are   unique.

Let us compare our results and those
from \cite{RZ_21} concerning the existence 
and uniqueness of strong solutions.
For simplicity we fix $T\in(0,\infty)$
and suppose that $b(t,x) =0$ for $t\not \in [0,T]$.

{\em Case $p>d$ (and $q>2$)\/}. Set $p_{b}=p_{D\sigma}=d$ ($\sigma$ is the unit matrix).
Then take a  constant $\hat N>0$
and let
$$
\lambda(t)=\hat N\Big(\int_{\bR^{d}}
|b (t,x)|^{p  }\,dx\Big)^{1/(p  -d )}.
$$
Also define  
$$
b_{M} (t, x)=b  (t, x)I_{|b  (t,x)|\geq \lambda(t)}
$$
and observe that  for $B\in\bB_{\rho}$  we have
$$
\dashint_{B }|b_{M} (t, x)|^{d}\,dx
\leq \lambda^{d-p  }(t)
\dashint_{B }|b  (t, x)|^{p  }\,dx
\leq N(d)\hat N^{d-p  }\rho^{-d}.
$$
 
Here $N(d)\hat N^{d-p  }$ is as small
as we like if   $\hat N$
large enough. Furthermore, for 
$b_{B} =b -b_{M} $ it holds that
 $|b_{B} |\leq \lambda $ and
$$
\int_{0}^{T}\lambda^{2}(t)\,dt
=\hat N^{2}\int_{0}^{T}\Big(\int_{\bR^{d}}|b (t,x)|^{p  }
\,dx\Big)^{q  /p  }\,dt<\infty.
$$
This shows that $b$ is admissible and $b_{M}$
is sufficiently small. Thus, ($\sigma$ is the unit matrix in \cite{RZ_21}) Theorem \ref{theorem 7.18.1}
is applicable and yields a strong solution
conditionally unique.

It turns out that there is a unique
(unconditional) solution if $p\geq d+1$. To show this
  set $\frp_{b}=p$, $\frq_{b}=q$,  choose
$p_{0}$ slightly less that $d$ $(=p_{b}$) and find
sufficiently large $q_{0}$ to satisfy \eqref{5.31.2}.
Then
 $$
 \frac{d}{p_{0}}+\frac{2}{q_{0}}<2,\quad
 \frac{d}{\frp_{b}}+\frac{2}{\frq_{b}}=1,\quad
\frac{d}{\frp_{b}}+\frac{1}{\frq_{b}}\leq 1.
$$  
Furthermore, \eqref{3.14.2} is satisfied with as small $\hat b$ as we like
 on the account of taking appropriate
$r_{b}$. Indeed,  
for any $C\in\bC_{r}$,
$$
\dashnorm b\|_{L_{\frp_{b},\frq_{b}}(C_{r})}
=Nr^{-d/\frp_{b}-2/\frq_{b}}\|b\|_{L_{\frp_{b},\frq_{b}}(C_{r})}
= Nr^{-1}\|b\|_{L_{p ,q }(C_{r})},
$$
where $N=N(d)$, and the last norm can be made   as small as we like on the account of choosing $r$ small enough. 
 
Thus, Theorem \ref{theorem 7.24.1} is applicable
(recall that $\sigma$ is the unit matrix).
In particular, if $p\geq d+1$, any solution is strong and unique and requiring property (b) is redundant.
Nevertheless, observe that by Theorem
\ref{theorem 6.1.1} property  (b) holds  albeit with some extra restrictions on $p,q$.

If $p\in(d,d+1)$, for a {\em strong\/} solution
$y_{\cdot}$ 
to coincide with the strong one $x_{\cdot}$, property (b), with $y_{\cdot}$ in place of $x_{\cdot}$,
is required in \cite{RZ_21}, where
$ p, q$ are {\em any\/} numbers satisfying \eqref{7.29.5}.
According to Theorem \ref{theorem 7.18.1}
we only need (a) with {\em some\/} $p,q$ satisfying \eqref{5.31.2} and do not need $y_{\cdot}$ to be a {\em strong\/} solution (it might be that this
discrepancy is caused by a sloppy use of the quantifiers). On the other hand, although the condition on the range of $p,q$
from \cite{RZ_21} is stronger (in uniqueness),   the right-hand
side of \eqref{7.29.4} is smaller than 
the right-hand side of \eqref{7.19.1} 
for $p=p_{0},q=q_{0}$.

{\em Case $p=d$ and $b\in C([0,T],L_{d})$\/}. In that case
\begin{equation}
                                  \label{6.3.1}
\lim_{r\downarrow 0}\sup_{t\in[0,T]}
\sup_{B\in \bB_{r}}\|b(t,\cdot) \|_{L_{d}(B)}=0.
\end{equation} 

Then, one can take $p_{b} =d$, choose any
$  p_{0} ,q_{0}  $ satisfying \eqref{5.31.2} and set $b_{M}=b $
which guarantees the arbitrary smallness of $\hat b_{M}$ for small 
enough $\rho_{b}$.    Therefore, this case is also covered by 
Theorem \ref{theorem 7.18.1} and the implications
of that are the same as in the  case of $p\in(d,d+1)$
(a result comparable
to the one in \cite{RZ_21}).

 We see that, actually, condition that
$b\in C([0,T],L_{d})$ can be replaced with $\leq
\varepsilon$ in
\eqref{6.3.1} in place of $=0$, for $\varepsilon>0$ small enough, which holds, for instance,
if the norms $\|b(t,\cdot)\|_{L_{d}}$ are uniformly
sufficiently small, that is imposed as one of alternative conditions
in \cite{RZ_21}.

Another, border case of \eqref{6.2.2}
when $p=\infty,q=2$, not covered in \cite{RZ_21},
also can be treated by our methods. We take
$b_{M}=0$, $p_{b}=d$,  
any $p_{0}\in (2\vee (d/2),d),q_{0}\in(1,\infty)$,
satisfying \eqref{5.31.2},  and observe that
\eqref{7.26.1} is satisfied with $f=b$.  In that case
Theorem \ref{theorem 7.18.1} is applicable
and yields a strong solution  which is unique
in the set of all solutions possessing property 
($\text{a}'$).

However, it is necessary to say that our results
do not cover all spectrum of results \cite{RZ_21},
one of the main emphasis of which is proving
weak differentiability of strong solutions with respect to the initial starting point.
On the other hand, in contrast to
\cite{RZ_21}, our diffusion matrix is not constant
or even continuous.

It is also worth noting a recent paper \cite{KM_23} which largely extends
{\em some\/} of the results of \cite{RZ_21} and the current paper to
the case of form-bounded drift when the diffusion is constant.
However, not all   results in \cite{RZ_21} and the current paper are covered because, for instance, of the assumption
in \cite{KM_23} that our $b_{B}\equiv0$. 
\end{remark}

\begin{remark}
                            \label{remark 6.4.4}
In Example 2.3 of \cite{10} there is a function $b=b(x)$,
which is admissible with   $p_{b} =d-1$
but is such that $b\not\in L_{ p_{b}+\varepsilon,\loc}$,
no matter how small $\varepsilon>0$ is. Therefore,
it is way beyond condition \eqref{6.2.2}, which is called the ``critical'' condition, and shows that we are considering
a ``supercritical'' case of $b$.  

We know from the example in Remark 3.15 of \cite{13}  
that the function
$$
f(t,x)=I_{1>t>0,|x|<1}|x|^{-1}\Big(\frac{|x|}{\sqrt t}\Big)^{1/(d+1)}
$$
does not satisfy \eqref{6.2.2} no matter
what $p,q$ are and does not satisfy
\eqref{3.14.2} if $\frp\geq\frq$ with $L_{p,q}$-norm
understood as in \eqref{6.4.1} but does satisfy
\eqref{3.14.2} if the norm  is taken from
\eqref{7.29.6} with any $\frp_{b},\frq_{b}$, such  that
$d/\frp_{b}+1/\frq_{b}=1$ and $d+1<\frq_{b}<2(d+1)$
implying that $\frp_{b}<\frq_{b}$,  and, as stated in Remark 3.15 of \cite{13}, the equation with unit $\sigma$
and $|b|=\varepsilon f$, for sufficiently
small $\varepsilon$, has a weak solution
starting from any point and this solution is weakly unique. 

We modify this example to have strong
solution which is  unique among all solutions, 
 owing to Theorem \ref{theorem 7.24.1}, in the following way. For sufficiently large integer  $n_{0}$ and $n\geq n_{0}$ define
$$
a_{n}=n^{-1}\ln^{-2}n,\quad b_{n}=(1+a_{n})a_{n}.
$$
By sufficiently large $n_{0}$ we mean that
$b_{n}\leq a_{n-1}\leq 1$ for $n\geq n_{0}$. Then set
$$
\kappa(t)=\sum_{n\geq n_{0}}I_{(a_{n},b_{n})}(t),\quad
g=\kappa f.
$$
It is not hard to check that $g$ still does not satisfy
\eqref{6.2.2} with $L_{p,q}$-norm
understood as in \eqref{6.4.1}.

Observe that the above $\frp_{b}$ satisfies
$2\vee(d/2)<\frp_{b}$ and, therefore, it is not hard to see that
there exist  $ p_{0}$, $q_{0}$ such that
$p_{0}<\frp_{b}$,
$q_{0}<\frq_{b}$ and \eqref{5.31.2} holds.  Now Theorem \ref{theorem 7.24.1} implies that the equation $dx_{t}=dw_{t}
+b(t,x_{t})\,dt$ with, say zero initial condition
and $|b|=\varepsilon g$ has a unique solution
and it is a strong one, provided that
$\varepsilon>0$ is small enough and $g$
is admissible with    any $p_{g}<d$.

 We  set $g_{M}(t,x)
=g(t,x)I_{|x|\leq\sqrt t}$ and  
   take any $p_{g}\in (1,d)$.  Since $g_{M}\leq |x|^{-1}$
and $p_{g}<d$,
for any $r>0$,
$$
\int_{B_{r}}|g_{M}(t,x)| ^{p_{g}}\,dx\leq 
\int_{B_{r}}|x| ^{-p_{g}}\,dx=Nr^{d-p_{g}}.
$$
If $|x|\leq 2r$, then
$$
\int_{B_{r}(x)}|g_{M}(t,y)| ^{p_{g}}\,dy\leq 
\int_{B_{3r} }|g_{M}(t,y)| ^{p_{g}}\,dy\leq Nr^{d-p_{g}}.
$$
In case  $|x|\geq 2r$, we have $g_{M}\leq r^{-1}$
on $B_{r}(x)$ and we get
$$
\int_{B_{r}(x)}|g_{M}(t,y)| ^{p_{g}}\,dy
\leq Nr^{d-p_{g}}
$$
again. This is part of conditions
needed for $ g$ to satisfy to qualify for being
an admissible function. The remaining part
is satisfied as well, since
$$
\sup_{x}g_{B}(t,x)\leq \kappa(t)t^{-1/2}
\sup_{|x|>\sqrt t}\big(\sqrt t/|x|\big)^{d/(d+1)}
=\kappa(t)t^{-1/2}
$$
and
$$
\int_{0}^{1}\kappa(t)t^{-1 }\,dt=
\sum_{n\geq n_{0}}\ln(b_{n}/a_{n})=
\sum_{n\geq n_{0}}\ln(1+a_{n})\leq \sum_{n\geq n_{0}} a_{n}<\infty.
$$
  This example shows that using somewhat
unnatural norms as in \eqref{7.29.6} could
be quite essential.

We see that the Ladyzhenskaya-Prody-Serrin condition
\eqref{6.2.2} is rather rough in what concerns the
existence of strong solutions for equations with singular drift.
 Also note that this example is not covered by
the results of \cite{KM_23} because their Assumption (A)
is not satisfied even with its extension described in Remark 2
of \cite{KM_23}. 

\end{remark}

\begin{remark}
                \label{remark 6.3.1}
Our results contain those in \cite{1} 
where the coefficients are independent of time.
Moreover, observe that if $b$ is admissible, 
$d/2<p_{b}\leq d$, 
and $b_{B}=0$, then
Assumption \ref{assumption 3.14.1} is automatically satisfied with $\frp_{b}=p_{b}$ and any
sufficiently large $\frq_{b}$ (satisfying \eqref{3.14.1}
and allowing to find $p_{0},q_{0}$)
since for  $C=[s,s+r^{2})\times B\in \bC_{r}$,
$r\leq r_{b}$,
we have
$$
r\dashnorm b \|_{L_{\frp_{b},\frq_{b}}(C)}\leq r
\sup_{t\in[0,T]}\dashnorm b(t,\cdot)\|_{L_{p_{b}}(B)}
\leq \hat b_{M}.
$$

\end{remark}
\begin{remark}
                            \label{remark 7.21.1}
We will be extensively using the results from 
\cite{10} and \cite{13}, where instead of assumptions
on $\sigma$ there are assumptions on $a=\sigma\sigma^{*}$
quite different form Assumption \ref{assumption 12.19.1}. More precisely in  \cite{10} and \cite{13} combined 
it is assumed that for $\rho\leq \rho_{a}$  (with some $\rho_{a}\in(0,1]$)
\begin{equation}
                                 \label{6.3.01}
a^{\# }_{\rho}:= \sup_{C\in\bC_{\rho}}\dashint_{C}|a (t,x)-a_{C}(t)|
\,dx dt \leq \theta,
\end{equation}
where
$$
a_{C}(t)=\dashint_{C}a (t,x)\,dxds 
\quad (\text{note}\,\,t\,\,\text{and}\,\,ds)
$$
and with $\theta>0$ depending only on $d,\delta,p_{b},p_{0},q_{0},\frp_{b},\frq_{b}$.

Note that by Poincar\'e's inequality for $C\in\bC_{\rho}$
$$
\dashint_{C}|a (t,x)-a_{C}(t)|
\,dx dt \leq N(d)\rho\dashint_{C}|Da|\,dxdt
\leq N(d,d_{1},\delta)\rho\dashint_{C}|D\sigma|\,dxdt
$$
$$
\leq N(d,d_{1},\delta)\rho\Big(
\dashint_{C}|D\sigma_{M}|\,dxdt+\sqrt{\beta_{D\sigma}(\rho^{2}})\Big).
$$

We also use the results from \cite{Kr_306} 
 where we   require $Da$ to be admissible.
In Remark \ref{remark 7.21.2} we saw that
this requirement is satisfied.

This discussion shows that the results in
\cite{10}, \cite{13}, and \cite{Kr_306}
(with   $p_{0}$ in place of $p$ in \cite{Kr_306} where
  $p_{D\sigma}$ is also involved) are applicable
in what concerns the conditions on $a$
if $\widehat{D\sigma_{M}},\beta_{D\sigma}(\rho^{2}_{\sigma})\leq \varepsilon_{\sigma}$, where
$\varepsilon_{\sigma}>0$ depends only on $d,d_{1},\delta,p_{b}$, $p_{D\sigma}$, $p_{0},q_{0}$, $\frp_{b}$, $\frq_{b} $. Actually, we need one more restriction on $\widehat{D\sigma_{M}}$ coming from
our 
Lemma \ref{lemma 6.13.1}, but it is also in terms of
$d,\delta$ and we assume that it is taken care of
by $\varepsilon_{\sigma}$. One more 
implication of this lemma is that we need the results
  in
\cite{10}, \cite{13}, and \cite{Kr_306} to be true
with $\delta/4$ there in place of $\delta$. Again we suppose
that $\varepsilon_{\sigma}$ is chosen appropriately.

In what concerns the conditions on $b$, the combined restrictions from \cite{10}, \cite{13}, and \cite{Kr_306} are that
$\hat b_{M},\hat b\leq\varepsilon_{b}$, where
$\varepsilon_{b}>0$
 depends only on $d,\delta,p_{b}$, $p_{0},q_{0}$, $\frp_{b},  \frq_{b}$,   $ \rho_{D\sigma}$.

\end{remark}

In order to use freely the results from
\cite{10}, \cite{13}, and \cite{Kr_306} 
 for $b$ and $\cSigma=\sqrt{a}$ in place of $\sigma$, we impose the following.

\begin{assumption}
                         \label{assumption 7.21.2}
We have $\widehat{D\sigma_{M}},\beta_{D\sigma}(\rho^{2}_{\sigma})\leq \varepsilon_{\sigma}$
and $\hat b_{M},\hat b\leq\varepsilon_{b}$,
where $\varepsilon_{\sigma},\varepsilon_{b}$
are introduced in Remark \ref{remark 7.21.1}.

\end{assumption}

Note that $r_{b}$ is {\em not\/} entering the arguments of 
$\varepsilon_{\sigma},\varepsilon_{b}$. Also, actually not in
every result below we need $\varepsilon_{\sigma}$ and
$\varepsilon_{b}$ to depend on all their listed arguments,
but not to overburden the exposition with not so important details
we decided against such specification. 

\mysection{Auxiliary results}
                            \label{section 6.4.1}

It is convenient in this section to replace
equation \eqref{6.15.2}  with
\begin{equation}
                                         \label{5.28.2}
 x_{s}=x +\int_{0}^{s}\cSigma (t+u,x_{u})\,dw _{u}+
\int_{0}^{s}b(t+u,x_{u})\,du,
\end{equation}
where  $w_{t}$ is a $d$-dimensional Wiener
process. We will be dealing with weak solutions
of \eqref{5.28.2}. Observe that each solution 
of \eqref{6.15.2} is also a solution of \eqref{5.28.2}
with $\cSigma =\sqrt a$ and the new $dw_{u}$ defined as $ \cSigma ^{-1}\sigma \,
dw _{u}$.  An advantage of considering \eqref{5.28.2}
in place of \eqref{6.15.2} is 
that mollifications of $\cSigma$ preserve its properties. It is also worth noting the well-known fact that any solution of \eqref{5.28.2} satisfies \eqref{6.15.2}
on an extended probability space (with different 
Wiener process).

We suppose that Assumptions \ref{assumption 12.19.1},
\ref{assumption 3.14.1}, and 
\ref{assumption 7.21.2}  are satisfied.
 
The constants $N$ below depend only on $d,\delta$,
$p_{b},p_{D\sigma},\frp_{b},\frq_{b}$, $\rho_{D\sigma} ,
 \rho_{b}$,  $p_{0},q_{0}$,    the
 values at $\infty$  of $\beta_{D\sigma}$,
$\beta_{b}$ and their moduli of continuity. 
When we write $N=N(T,...)$, the constant is also allowed to depend on $T,...$. 
By $\hat N,\hat N(T,...)$ we denote the constants
which are also allowed to depend on $r_{b}$.

Define ($a=\cSigma ^{2}$) 
$$
\cL u(t,x)=\partial_{t}u(t,x)+(1/2)a^{ij}(t,x)
D_{ij}u(t,x)+b^{i}(t,x)D_{i}u(t,x),
$$
$$
\bR^{d}_{T}=[0,T]\times \bR^{d}.
$$
 
Here the main emphasis is on constructing
an  evolution family $T_{s,t}$ of operators
corresponding to \eqref{5.28.2}. 

\begin{lemma}
                       \label{lemma 4.28.1}
Additionally assume that 
$b$ and $Da$ are   bounded.
Then
for any $T\in(0,\infty)$ and bounded smooth $f$
on $\bR^{d } $ with compact support, there exists a unique $W^{1,2}_{p_{0}}(\bR^{d}_{T})$-solution  $u$ of
$\cL u=0$ in $\bR^{d}_{T}$ such that $u(T,\cdot)=f$.
Furthermore, for this solution and $t<T$
\begin{equation}
                     \label{4.28.1}
 \|u(t,\cdot)\|_{L_{ p_{0}}}\leq N(T,p)\|f\|_{L_{ p_{0}}},\quad \|Du(t,\cdot)\|_{L_{p_{0}}}
\leq N(T,p)(T-t)^{1/ p_{0}-1}\|f\|_{L_{ p_{0}}}, 
\end{equation}
\begin{equation}
                             \label{4.27.3}
|u(t,x)|\leq N(T)(T-t )^{1/q_{0}-1/\alpha }
\int_{\bR^{d}}\|f\|_{L_{ p_{0}}(B_{1}(z))}
\Psi(T-t,z-x)\,dz,
\end{equation}
where   $\alpha=\alpha(p_{0},q_{0})>1$
and, with $N=N(T)$,
$$
\Psi(t,x )=I_{|x|<2}+e^{-|x| /(Nt)}.
$$

\end{lemma}

Proof. The existence and uniqueness is a classical result.   Then,
the first estimate in \eqref{4.28.1} follows
directly from Theorem 1.7 of \cite{Kr_306}, the possibility to rewrite $\cL$ in the divergence form,
and the fact that $p_{b}>p_{0}\geq2$ (in fact, $p_{0}>2$). Formally speaking we should have added that we also  needed
$p_{0}<   p_{D\cSigma}$, which we did not require.  However,
at this point this is irrelevant because once $u$ is smooth one can write the equation for $|u|^{p}$
and then proceed estimating $\|u(t,\cdot)\|_{L_{ p_{0}}}$
as in \cite{Kr_306}. This comment is also relevant to
the estimate of $\|Du(t,\cdot)\|_{L_{p_{0}}}$.

To prove it,
take an $s\in(0,T)$, set $s_{0}=(s+T)/2$, and take
 a smooth decreasing function $\zeta(t)$, such that
 $\zeta(t)=1$ for $t\leq s $,
$\zeta(t)=0$ for $t\geq s_{0}$, and set $v=u\zeta $. Then $\cL v-u\zeta'=0$  and we use Theorem 1.11 of \cite{Kr_306} with
one alteration that, in the definition of $K$ in its statement,
$Df_{s}$ could be dropped because the estimate
of the term $I_{8}$ in its proof could be done as follows:
$$
I_{8}:=\theta_{t}|v_{t}|^{p-1}f_{t(\eta)}+
N|\eta|^{2}|v_{t}|^{p-2}|Dg_{t}|^{2}\leq \theta_{t}|v_{t}|^{p-1}f_{t(\eta)}+ V_{t}^{2}+N|\eta|^{p}|Dg_{t}|^{p}, 
$$
where
$$
\theta_{t}|v_{t}|^{p-1}f_{t(\eta)}\sim_{\kappa}-(p-1)v_{t}|^{p-2} v_{t(\eta) }
f_{t}
$$
$$
\leq \varepsilon \big(|v_{t}|^{(p-2)p/(2p-2)}\big)\Big(
|v_{t}|^{(p-2)p/(2p-2)}\big(|\eta| |D^{2}u|\big)^{p/(p-1)}\Big)+\varepsilon^{-1}|\eta|^{p}
|f_{t}|^{p}
$$
$$
\leq \varepsilon|v_{t}|^{p}+\varepsilon|\eta|^{p}|Du_{t}|^{p-2}
|D^{2}u|^{2}+\varepsilon^{-1}|\eta|^{p}
|f_{t}|^{p}
$$
Then we find  
$$
 \int_{\bR^{d}}|Dv(t,x)|^{ p_{0}}\,dx
\leq N\int_{t}^{T}\int_{\bR^{d}}|u\zeta'(t,x)|^{ p_{0}}\,dx
dt
$$
$$
\leq N\sup|\zeta'|^{ p_{0}-1}\int_{t}^{T}\int_{\bR^{d}}|u|^{ p_{0}}|\zeta'|\,dx
dt
\leq N(T-s)^{-( p_{0}-1)}\|f\|_{L_{ p_{0}}}^{ p_{0}},
$$
where the last inequality follows from the first estimate
in \eqref{4.28.1} and the fact that the integral
of $|\zeta'|$ is one. This yields the second
estimate in \eqref{4.28.1} for $t=s$,
since $v(s,x)=u(s,x)$.

Next,  there is a Markov diffusion process
corresponding to $\cL$ (see \cite{10}) in terms of which
$$
v(t,x)=E_{t,x}\int_{0}^{T-t}u\zeta'(t+s,x_{s})\,ds.
$$
Therefore, 
it follows from   Theorem 3.9 of \cite{10}
that (here we use that   $1<p_{0}<p_{b},d/ p_{0}+2/q_{0}<2$) we have
$$
|v|\leq N
\Big(\int_{0}^{T}\Big(\int_{\bR^{d}}|u|^{ p_{0}}\,dx\Big)^{q_{0}/ p_{0}}|\zeta'|^{q_{0}}\,dt\Big)^{1/q_{0}}
\leq N\|f\|_{L_{ p_{0}}} \Big(\int_{0}^{T}|\zeta'|^{ q_{0}}\,dt\Big)^{1/ q_{0}}
$$
$$
\leq N\|f\|_{L_{ p_{0}}}(T-s)^{1/ q_{0}-1}.
$$
This yields
\begin{equation}
                               \label{4.28.4}
|u(s,x)|\leq N(T-s )^{1/ q_{0}-1}
 \|f\|_{L_{ p_{0}} }.
\end{equation}
 
On the other hand, take a nonnegative $\eta\in C^{\infty}_{0}(B_{1})$ with unit integral 
and set $\eta_{z}(x)=\eta(x-z)$. Then
$$
u(t,0)=\int_{\bR^{d}}E_{t,0}f\eta_{z}(x_{T-t})\,dz,
$$
\begin{equation}
                               \label{4.28.3}
|u(t,0)|\leq N\int_{\bR^{d}}(E_{t,0}|f|^{\alpha}\eta_{z}(x_{T-t}) )^{1/\alpha}P^{1/\alpha'}_{t,0}(x_{T-t}\in B_{1}(z))\,dz,
\end{equation}
where $\alpha'=\alpha/(\alpha-1)$ and  $\alpha=\alpha(p_{0},q_{0})>1$
is such that $d/p_{0}+2/q_{0}<2/\alpha$ and  
$1<p_{0}/\alpha< p_{b} $. 
By applying \eqref{4.28.4} to $p_{0}/\alpha$ 
and $q_{0}/\alpha$ in place of $p_{0}$ and $q_{0}$,   we find
$$
(E_{t,0}|f|^{\alpha}\eta_{z}(x_{T-t}) )^{1/\alpha}\leq
N (T-t)^{1/q_{0}-1/\alpha}
 \|f\|_{L_{p_{0}}(B(z))}.
$$

Furthermore,
as it is pointed out in \cite{10}, all results of 
\cite{Kr_pta} are applicable
(with $\rho_{b}$ there defined in Theorem 3.1 of 
\cite{10}
as a function of (our) $\rho_{b},d,\delta$, and the modulus
of continuity of $\beta_{b}$). In particular,
 Theorem 3.8 of \cite{Kr_pta} implies that for
$|z|\geq2$
$$
P _{t,0}(x_{T-t}\in B_{1}(z))\leq Ne^{-|z| /(NT-Nt)}.
$$
This and \eqref{4.28.3}  prove \eqref{4.27.3}
at $x=0$, which is certainly enough.
The lemma is proved. \qed

Next, recall that $a=\sigma\sigma^{*}$, where $\sigma$ is taken from Section
\ref{section 12.20.1},  and  take the
Markov diffusion process $X=((t+s ,x_{s}),
\infty,\cN_{s},P_{t,x})$ constructed in 
Theorems 3.9 and 4.6 of \cite{10}
 by approximating
$a, b$   by mollified functions with mollification acting only in $x$ (first method).
The trajectories $(t+s,x_{s})$ of the process $X$ are constructed
as solutions of \eqref{5.28.2} with $\cSigma=\sqrt a$.

\begin{remark}
                           \label{remark 7.24.1}

By Theorem 3.9 of \cite{10}
 property (a) from Theorem \ref{theorem 7.18.1}
holds
with $E=E_{0,0}$ and $N=N(T)$ (depending also on other data as 
agreed upon above).
\end{remark}

For the process $X$ define
$$
T_{s,t}f( x)=E_{s,x}f( x_{t-s}).
$$
Then, by using the passage to the limit from
regular coefficients to the original $a$ and $b$
and using Lemma \ref{lemma 4.28.1} we obtain the following.

\begin{lemma}
                       \label{lemma 4.30.1}
Take $t_{0}\in(0,\infty)$, $f\in L_{ p_{0}}$, and set
$u(t,x)=T_{t,t_{0}}f( x)$, $t\in[0,t_{0})$. Then, for each $t$
the Sobolev derivative $Du(t,x)$ exists and
\eqref{4.28.1} and \eqref{4.27.3} hold.

\end{lemma}

The main step in proving Theorem
\ref{theorem 6.1.1} is done in the following.  
\begin{lemma}
                             \label{lemma 6.1.2}
For any $T\in(0,\infty)$ and $p,q$ satisfying
\eqref{6.1.2} and  Borel nonnegative $f$ we have
$$
E_{0,0} \int_{0}^{T}f( s,x_{s})\,ds
$$
\begin{equation}
                                \label{6.1.4}
\leq 
N(T,p,q)\Big(\int_{0}^{T}\Big(\int_{\bR^{d}} \|f(s,\cdot)\|_{L_{p }(B_{1}(z))}^{p}\Psi(s,z )\,dz\Big)^{q/p}ds\Big)^{1/q},
\end{equation}
where $\Psi(t,x)=I_{B_{1}}+e^{-|x| /(Nt)}$.
\end{lemma}

Proof. Set
  $\kappa=p/p_{0}\geq1$. Then we use Lemma
\ref{lemma 4.30.1}  
to see that for each $s\in (0,T] $
$$
\big(E_{0,0}f(s,x_{s})\big)^{\kappa}\leq 
E_{0,0}f^{\kappa}(s,x_{s})\leq 
Ns^{1/q_{0}-1/\alpha}\int_{\bR^{d}}
\|f^{\kappa}(s,\cdot)\|_{L_{p_{0}}(B_{1}(z))}\Psi(s,z)\,dz,
$$
where by H\"older's inequality the integral is dominated by
$$
N\Big(\int_{\bR^{d}}
\|f^{\kappa}\|^{p_{0}}_{L_{p_{0}}(B_{1}(z))}\Psi(s,z)\,dz
\Big)^{1/p_{0}}=N\Big(\int_{\bR^{d}}
\|f(s,\cdot) \|^{p}_{L_{p}(B_{1}(z))}\Psi(s,z)\,dz
\Big)^{\kappa/p}
$$
since
 $\int_{\bR^{d}}
 \Psi(s,z)\,dz\leq N$. It follows that
$$
E_{0,0} \int_{0}^{T}f( s,x_{s})\,ds
\leq N\int_{0}^{T}s^{(1/q_{0}-1/\alpha)/\kappa}
\Big(\int_{\bR^{d}}
\|f(s,\cdot) \|^{p}_{L_{p}(B_{1}(z))}\Psi(s,z)\,dz
\Big)^{1/p}\,ds.
$$
Then using again H\"older's inequality after
observing that $(q/(q-1))(1/\alpha-1/q_{0})/\kappa<1$
for $p,q$ satisfying  \eqref{6.1.2}.
The lemma is proved. \qed

We now need to introduce new parameters
   $\beta_{0},\beta_{0}' $ into the play. We choose
    $\beta_{0}\in(1,2)$ so that
\begin{equation}
                           \label{3.14.02}
\beta_{0}<\frp_{b}\wedge \frq_{b},\quad \beta_{0}p_{0}\leq\frp_{b},\quad
\beta_{0}q_{0}\leq \frq_{b},\quad
2(\beta_{0}-1) p_{0}<
\frp_{b}\wedge \frq_{b}, 
\end{equation}
and set
 \begin{equation}
                           \label{7.24.2} 
\sfp=\frp_{b}/\beta_{0},\quad\sfq=\frq_{b}/\beta_{0}.
\end{equation}

Observe that
\begin{equation}
                           \label{7.26.01}
\frac{d}{\sfp}+\frac{2}{\sfq }\geq\beta_{0}>1. 
\end{equation}
The same conditions we impose on $\beta_{0}'$ with
the addition that $\beta_{0}'<\beta$. We concider
$\beta_{0}$. $\beta'_{0}$ as functions of $p_{0},q_{0},
\frp_{b},\frq_{b}$, and this is the only reason they are not included
in the set of arguments of $\varepsilon_{\sigma}.\varepsilon_{b}$.

Owing to Theorem 3.9
of \cite{10} first applied to approximating processes and the passing to the limit  for any $T\in(0,\infty)$,
$m=1,2,...$, and Borel nonnegative $f$ on $\bR^{d+1}$
and $(t,x)$ 
\begin{equation}
                          \label{7.19.3}
E_{t,x}\Big(\int_{0}^{T}f(s,x_{s})\,ds\Big)^{m}
\leq N(T,m)\|f\|^{m}_{L_{p_{0},q_{0}}}.
\end{equation}
This and the fact that $p_{0}\leq \sfp$,
$q_{0}\leq \sfq$ (see \eqref{3.14.02}, \eqref{7.24.2}) imply that for any $R\in(0,\infty)$
and Borel nonnegative $f$ on $\bR^{d+1}$
\begin{equation}
                             \label{12.11.06}
E_{t,x}\int_{0}^{\tau_{R}}f(s,x_{s})\,ds\leq N(T,R)\|fI_{C_{R}(t,x)}\|_{L_{\sfp ,\sfq}}
\leq N(T,R)\|f \|_{E_{\sfp ,\sfq,\beta_{0}}},
\end{equation}
where $\tau_{R}$ is the first exit time of $(s,x_{s})$ from $C_{R}(t,x)$.

  Therefore, in light of condition \eqref{3.14.1}, by the construction in Theorem 3.18 and   weak uniqueness  Theorem 3.8
of \cite{13} (here we use, in particular, that   $\hat b\leq \varepsilon_{b}$), the same Markov process $X$
(in the sense of finite-dimensional distributions)
can be constructed by using the mollification of $a$ and $b$ in $(t,x)$ (second method).

The appearance of $\hat N$ in the following lemma
needs an explanation. In the constants entering the estimates 
in \cite{13} no $r_{b}$ appears. The reason for this
is that in \cite{13} the parabolic dilation is
used to reduce $r_{b}\leq 1$ to $1$. To apply the results
obtained after such dilation we need to make
the inverse contraction and this leads to the
appearance of $r_{b}$ in the constants' arguments.

\begin{lemma}
                       \label{lemma 4.30.2}
Take $t_{0}\in(0,\infty)$, $f\in C_{0}^{\infty}$, and set
$u(t,x)=T_{t,t_{0}}f( x)$, $t\in[0,t_{0})$. Then  
the Sobolev derivatives $Du, D^{2}u,\partial_{t}u$
exist, $\cL u=0$ for $t<t_{0}$, and
\begin{equation}
                         \label{4.30.04}
\|u,Du,D^{2}u,\partial_{t} u\|_{E_{\sfp,
\sfq,\beta_{0}}(\bR^{d}_{t_{0}})}
\leq \hat N(t_{0})\|f\|_{W^{2}_{  p}},
\end{equation}
where $  p=(d/\beta_{0})\vee \sfp$.

\end{lemma}

Proof. Modify $\cL$ after time $t_{0}$ by setting
$\cL=\partial_{t}+\Delta $ and
define $g=-I_{t>t_{0}}\cL (f\zeta)$, where $\zeta(t)$ is smooth $\zeta(t_{0})=1$ and $\zeta(t)=0$ for $t\geq 2t_{0}$.
Then take sufficiently large $\lambda>0$
and consider the equation
\begin{equation}
                         \label{4.30.05}
\lambda v-\cL v=e^{ \lambda t}g.
\end{equation}
Since, obviously, $e^{ \lambda t}g\in L_{\sfp,
\sfq,\beta_{0}}$, by Theorem 3.7 of \cite{13}
there is a unique $v\in E^{1,2}_{\sfp,
\sfq,\beta_{0}}$
satisfying \eqref{4.30.05} in $\bR^{d+1}$.
By uniqueness $e^{-\lambda t}v(t,\cdot)=f\zeta(t,\cdot)$ for
$t\geq t_{0}$ and $v(t_{0},\cdot)=fe^{\lambda t_{0}}$.
Observe that since $\beta<2$, $v$ is bounded  
(see Lemma 2.5 of \cite{11}).
By   It\^o's formula (Lemma 3.6  of \cite{13},
here one needs to apply the previous arguments to
$\beta' $   in place of $\beta$) for $t<t_{0}$
$$
v(t,x)=E_{t,x}e^{-\lambda (\tau_{R}\wedge (t_{0}-t))}
v( t+\tau_{R}\wedge (t_{0}-t),x_{ \tau_{R}\wedge (t_{0}-t)}),
$$
where $\tau_{R}$ is the first exit time of
$ x_{s} $ from $B_{R}(x)$. As $R\to\infty$,
$\tau_{R}\to\infty$, $ \tau_{R}\wedge (t_{0}-t)
\to t_{0}-t$, and by the dominated convergence theorem
$$
v(t,x)=e^{-\lambda(t_{0}-t)}Ev(t_{0},x_{t_{0}-t})=e^{\lambda t}u(t,x).
$$
To finish the proof, now it only remains to
use the estimate of the $E^{1,2}_{\sfp,
\sfq,\beta_{0}}$-norm of $v$ from Theorem
3.7 of \cite{13} and observe that  the $E _{\sfp,
\sfq,\beta_{0}}$-norm of $g$
is dominated by the $W^{2}_{ p}$-norm of $f$
in light of H\"older's inequality.
The lemma is proved.\qed
 
In the next lemma we allow $f$ to be more general.

\begin{lemma}
                       \label{lemma 5.1.1}
Let $t_{0}\in(0,\infty)$, $f\in L_{ p_{0},\loc}$,
and 
$$
\sup_{z\in\bR^{d}}\|fI_{B_{1}(z)}\|_{L_{ p_{0} }}<\infty,
$$ 
 and set
$u(t,x)=T_{t,t_{0}}f( x)$, $t\in[0,t_{0})$. Then  
the Sobolev derivatives $Du$, $D^{2}u$, $\partial_{t}u$
exist, $\cL u=0$ for $t<t_{0}$, and for any $t<t_{0}$  
\begin{equation}
                         \label{4.30.4}
\|u,Du,D^{2}u,\partial u\|_{L_{\sfp,
\sfq,\beta_{0}}(\bR^{d}_{t})}
\leq \hat N(t_{0} )(t_{0}-t)^{1/q_{0}-1/\alpha}
\sup_{z\in\bZ^{d}}\|fI_{B_{1}(z)}\|_{L_{ p_{0} }}.
\end{equation}

\end{lemma}

Proof. First take
$f\in C^{\infty}_{0}$, take
$v$ from the proof of Lemma \ref{lemma 4.28.1}
(with $T=t_{0}$),
note that $v\in E^{1,2}_{\sfp,\sfq,\beta_{0}}(\bR^{d}_{\infty})$ by Lemma \ref{lemma 4.30.2},
$\cL v-\lambda v=\zeta'u-\lambda v$ for any $\lambda$, which by Theorem 3.5 of \cite{11}
shows that the $E^{1,2} _{\sfp,
\sfq,\beta_{0}}(\bR^{d}_{\infty})$-norm of $v$
is dominated by the sum of the $E_{\sfp,
\sfq,\beta_{0}}(\bR^{d}_{\infty})$-norms of $v$ and $\zeta'u$, which
is, obviously less than $\sup|\zeta u|+\sup|\zeta'u|$, that
admits an estimate like \eqref{4.27.3}.
The case of general $f $
is treated by approximation.
The lemma is proved. \qed

 \begin{theorem}
                         \label{theorem 5.15.1}
For the function $u$ from Lemma \ref{lemma 5.1.1},
for any $r <\sqrt t_{0}$ and $x\in\bR^{d}$,  
$$
\| u,Du,D^{2}u,\partial_{t} u\|_{L_{\sfp,
\sfq,\beta_{0}}(C_{r }(0,x))}
$$
\begin{equation}
                         \label{5.15.1}
\leq \hat N(t_{0})(t_{0}-r^{2} )^{-2+1/q_{0}-1/\alpha}
\int_{\bR^{d}}\|f\|_{L_{ p_{0}}(B_{1}(z ))}
\Psi(t_{0}-r^{2},z-x )\,dz.
\end{equation}
\end{theorem}

Proof. Set $r^{2}_{0}=(r^{2}+t_{0})/2$, $r_{0}>0$.
By a standard way of localizing the estimate
in Theorem 3.5 of \cite{11} (see, for instance, 
Lemma 2.4.4 of \cite{Kr_08}) using the interpolation
inequalities from Lemma 5.10 of \cite{Kr_22} one arrives at the following
\begin{equation}
                          \label{5.15.3}
\|D^{2}u,\partial_{t} u\|_{L_{\sfp,
\sfq,\beta_{0}}(C_{r} )}\leq \hat N(t_{0})(t_{0}-r^{2})^{-2}\sup_{C_{r_{0}} }|u|.
\end{equation}
By using Lemma \ref{lemma 4.30.1}, we come to
\eqref{5.15.1} with $x=0$ with no $u$ and $Du$.
The estimate of $u$ is obtained from \eqref{4.27.3}
and $Du$  is taken care of by interpolation inequalities.
Then we have \eqref{5.15.1} as is with $x=0$
which is enough and the theorem is proved.\qed

\begin{theorem}
                             \label{theorem 5.29.1}
Under the assumption  of Lemma \ref{lemma 5.1.1}
for any $r<\sqrt t_{0}$
\begin{equation}
                          \label{5.29.2}
\|Du\|_{L_{2 p_{0}}(\bR^{d}_{r})}\leq \hat N(t_{0},r)
\|f\|_{L_{ p_{0}}}.
\end{equation}
\end{theorem}

Proof. Let $(p,q)(\beta_{0}-1)=(\sfp,\sfq)\beta_{0}$ 
($=(\frp_{b},\frq_{b} )$). Owing to interpolation inequalities
 (Theorem 5.11 of \cite{Kr_22}, here we use \eqref{7.26.01}) 
\eqref{5.15.1} implies that the norm of $Du$ in $E_{p,q,\beta_{0}-1}
(C_{r}(0,x))$ is dominated by the right-hand side
of \eqref{5.15.1}. In light of \eqref{3.14.02}
we have $2(\beta_{0}-1) p_{0}<
\frp_{b}\wedge \frq_{b}$ implying that   
$p,q\geq2 p_{0}$. It follows by H\"older's
inequality that
\begin{equation}
                           \label{7.19.5}
\|Du\|^{2 p_{0}}_{L_{2 p_{0}}(C_{r }(0,x)))}
\leq \hat N(t_{0},r)\Big(
\int_{\bR^{d}}\|f\|_{L_{ p_{0}}(B_{1}(z ))}
\Psi(t_{0}-r^{2},z-x )\,dz\Big)^{2 p_{0}}.
\end{equation}
Here the integral over $x\in \bR^{d}$ of the left-hand sides is dominated the $2 p_{0}$'th power of
the left-hand side of \eqref{5.29.2} and
the integral of the right-hand sides is less than
a constant $\hat N(t_{0},r)$ times
$$
\int_{\bR^{d}}\int_{\bR^{d}}\|f\|^{2 p_{0}}_{L_{ p_{0}}(B_{1}(z ))}
\Psi(t_{0}-r^{2},z-x ) \,dzdx\leq N(t_{0},r)
\int_{\bR^{d}}\|f\|^{2 p_{0}}_{L_{ p_{0}}(B_{1}(z ))}\,dz
$$
$$
\leq N(t_{0},r)\|f\|^{  p_{0}}_{L_{ p_{0}} }
\int_{\bR^{d}}\|f\|^{  p_{0}}_{L_{ p_{0}}(B_{1}(z ))}\,dz
\leq N(t_{0},r)\|f\|^{2  p_{0}}_{L_{ p_{0}} }.
$$
This proves the theorem.\qed

\begin{corollary}
                        \label{corollary 7.19.1}
If $R\in(0,\infty)$ and $f=0$ outside $B_{R}$, then
\begin{equation}
                           \label{7.19.4}
\int_{0}^{r^{2}}\|Du(t,\cdot)I_{B^{c}_{2(R+1)}}\|_{L_{2 p_{0}}}^{2 p_{0}}\,dt\leq \hat N(t_{0},r)e^{-R/N}\|f\|^{2 p_{0}}_{L_{ p_{0}}}.
\end{equation}
 
\end{corollary}

Indeed, the left-hand side is less than a
constant, depending only on $r$, times the integral
of the left-hand side of \eqref{7.19.5} over $B_{2(R+1)}^{c}$, which is less than $\hat N(t_{0},r)$ times
$$
\int_{|x|\geq 2(R+1)}
\Big(
\int_{\bR^{d}}\|f\|_{L_{ p_{0}}(B_{1}(z ))}
\Psi(t_{0}-r^{2},z-x )\,dz\Big)^{2 p_{0}}\,dx
$$
$$
\leq\|f\|^{2 p_{0}}_{L_{ p_{0}}}\int_{|x|\geq 2(R+1)}
\Big(
\int_{|z|\leq R+1} 
\Psi(t_{0}-r^{2},z-x )\,dz\Big)^{2 p_{0}}\,dx
$$
Now \eqref{7.19.4} follows after noting that
in the integrals $|z-x|\ge (1/2)|x|$ and the
interior integral is less than a constant
of type $N$ in \eqref{7.19.4} times
$(R+1)^{d}e^{-|x|/N}\leq e^{-|x|/N}$ (with a different $N$).

Next, we consider a sequence of $a^{n},b^{n}$,
$n=1,2,...$,
having the same meaning as $a,b$, satisfying
the same assumptions as $a,b$ with the same
defining objects (in particular, $\delta$, $p_{a}$, $p_{b}$, $\rho_{Da}$,
$\rho_{b}$, the moduli
of continuity of $\beta_{Da},\beta_{b}$, and their values at $\infty$)
and satisfying the following,
which is supposed to hold throughout the rest of the section. Set $\sigma^{n}=\sqrt{a^{n}}$.

\begin{assumption}
                    \label{assumption 5.16.1}
For each $n$, the functions $b^{n},Da^{n}$ are bounded, $\sigma^{n},b^{n}$ satisfy Assumption \ref{assumption 7.21.2}, $a^{n}\to a$, $b^{n}_{B}\to b_{B} $ (a.e.)
as $n\to \infty$ and for each $ r>0$, $t\in\bR$
\begin{equation}
                               \label{5.16.2}
\int_{0}^{r}\|b^{n}_{M}(t+s,\cdot)-b_{M}(t+s,\cdot)\|_{L_{ p_{b}}(B_{r})}\,ds\to0
\end{equation}
as $n\to \infty$.
\end{assumption}

It is probably worth mentioning that  $\|b^{n}_{M}(t ,\cdot)\|_{L_{ p_{b}}(B_{r})}$ are uniformly boun\-ded on $\bR$ for each $r$.

Let $X^{n}=((t+s,x _{s}),\infty,\cN _{s},P^{n}_{t,x})$ be Markov processes corresponding to
$(a^{n},b^{n})$.  Then exactly as in the proof of Theorems 3.9  of \cite{10} one shows that for each $(t,x)$ there is a subsequence $n'\to \infty$
such that the $P^{n'}_{t,x}$-distributions of $(t+\cdot,x _{\cdot})$ on $C([0,\infty),\bR^{d+1})$ converge weakly to the distribution of a solution
of \eqref{5.28.2}. Since for each $n$,
estimate
\eqref{7.19.3} holds for $X^{n}$ with a constant  independent of $n$, it also holds for the limiting process, and then its distribution is uniquely
defined in light of Theorem 4.4 of \cite{10}. It follows that for the whole sequence we have that the $P^{n }_{t,x}$-distributions of $(t+\cdot,x _{\cdot})$ on $C([0,\infty),\bR^{d+1})$ converge weakly to the $P _{t,x}$-distribution
of $(t+\cdot,x _{\cdot})$. In particular,
(using $E^{n}$ as the expectation sign for $X^{n}$)
for $T^{n}_{s,t}f( x):=E^{n}_{s,x}f(x_{t-s})$, $0\leq
s\leq t<\infty $,
and bounded continuous $f$ on $\bR^{d}$ we have 
\begin{equation}
                               \label{5.16.3}
T^{n}_{s,t}f(x)\to T _{s,t}f(x)=:T^{0}_{s,t}f(x)
\end{equation}
as $n\to \infty$ for any $x\in\bR^{d} $.

\begin{lemma}
                            \label{lemma 5.17.1}
 Let $f\in L_{ p_{0}}$, $t_{0}\in(0,\infty)$, $S\in(0,t_{0})$. Then the convergence in
\eqref{5.16.3} holds uniformly with respect to $(s,x)
\in \bR^{d}_{S}$. Furthermore,
uniformly on $[0,S)$
\begin{equation}
                               \label{5.17.2}
\|T^{n}_{s,t_{0}}f -T^{0}_{s,t_{0}}f \|_{L_{ p_{0}}}\to0
\end{equation}
as $n\to \infty$. Finally, if $f$ is bounded,
the convergence in \eqref{5.16.3} holds uniformly with respect to $(s,x)$ in any bounded subset of
$\bR^{d}_{S}$.
\end{lemma}
 
Proof. Set $u^{n}(s,x)=T^{n}_{s,t_{0}}f( x)$. By Lemma \ref{lemma 5.1.1}
the $E^{1,2}_{\sfp,\sfq,\beta_{0}}(\bR^{d}_{S})$-norms
of $u^{n} $
are uniformly bounded, which by Lemma 2.8
of \cite{11} implies that they are uniformly
H\"older continuous in $\bR^{d}_{S}$. Therefore,
the convergence in \eqref{5.16.3} holds uniformly
on any compact set in $\bR^{d}_{t_{0}}$.
In light of \eqref{4.27.3} for
$n\geq0$, $x\in\bR^{d}$, and $s\leq S$
\begin{equation}
                                \label{5.20.1}
|u^{n}(s,x)|^{ p_{0}}\leq N(t_{0},S)\int_{\bR^{d}}
\|f\|^{ p_{0}}_{L_{ p_{0}}(B_{1}(z))}\Psi(t_{0}-S,z-x)\,dz.
\end{equation}
Here
$$
\int_{\bR^{d}}
\|f\|^{ p_{0}}_{L_{ p_{0}}(B_{1}(z))}\,dz<\infty
$$
and $\Psi(t_{0}-S,z-x)\to0$ as $|x|\to\infty$, so that
$u^{n}(s,x)\to0$   uniformly as $|x|
\to\infty$ and this yields their uniform
convergence to $u(s,x)$ in $\bR^{d}_{S}$.

To prove \eqref{5.17.2} again use \eqref{5.20.1}
to see that
for any $R>0$, ,
$$
\|u^{n}(s,\cdot)I_{B^{c}_{R}}\|^{ p_{0}}_{L_{ p_{0}}}
\leq N\int_{\bR^{d}}
\|f\|^{ p_{0}}_{L_{ p_{0}}(B_{1}(z))}\int_{|x|>R}\Psi(t_{0}-S,z-x)\,dxdz.
$$
This and the uniform convergence of $u^{n}$ prove 
\eqref{5.17.2}. 

To prove the last assertion of the lemma, take $R>0$
and set $g_{R}=fI_{B_{R}}$, $h_{R}=fI_{B_{R}^{c}}$,
$$
v^{n}_{R}(s,x)=T^{n}_{s,t_{0}}g_{R}(x),\quad w
^{n}_{R}(s,x)=T^{n}_{s,t_{0}}h_{R}(x). 
$$
 Then the first assertion of the lemma implies the said convergence
of $v^{n}_{R}$ to $v^{0} _{R}$. 
After that in only remains to observe that,
thanks to \eqref{4.27.3}, in any bounded
subset of $\bR^{d}_{S}$ the functions $w^{n}_{R}$
are uniformly small if $R$ is large.
The lemma is proved. \qed

\begin{theorem}
                           \label{theorem 5.22.1}
Let $f\in L_{ p_{0}}$, $t_{0}\in(0,\infty)$, $S\in(0,t_{0})$, then 
in the notation from the proof of Lemma \ref{lemma 5.17.1}
\begin{equation}
                       \label{5.22.1}
\|Du^{n}-Du^{0} \|_{L_{ p_{0}}(\bR^{d}_{S})}\to0
\end{equation}
as $n\to\infty$.
\end{theorem}

Proof. In light of \eqref{4.28.1} we may assume
that $f$ has compact support. Corollary 
\ref{corollary 7.19.1} shows that it suffices to prove that
\begin{equation}
                       \label{7.19.6}
\|Du^{n}-Du^{0} \|_{L_{ p_{0}}((0,S)\times B_{R} )}\to0
\end{equation}
as $n\to\infty$ for any finite $R$. Then
in light of Theorem \ref{theorem 5.29.1} it suffices
to show that $Du^{n}-Du^{0}\to 0$ in measure
as $n\to\infty$.

By embedding theorems 
$$
\| I_{B_{R}}(Du^{n}- Du^{0}) \|_{L_{\sfp,\sfq}(\bR^{d}_{S})}
$$
$$
\leq 
\varepsilon \| I_{B_{R}}(D^{2}u^{n} ,  D^{2}u^{0} ,\partial_{t}u^{n},\partial_{t}u^{0})\, \|_{L_{\sfp,\sfq}(\bR^{d}_{S})} +  N(\varepsilon,S,R)
\| I_{B_{R}}( u^{n}-  u^{0}) \|_{L_{\sfp,\sfq}(\bR^{d}_{S})} ,
$$
where $\varepsilon\in(0,1)$ is arbitrary.
Here
the factor of $\varepsilon$ is bounded by a constant independent of $n$ due to Theorem \ref{theorem 5.15.1}
and the last term tends to zero
as $n\to \infty$ owing to Lemma
\ref{lemma 5.17.1}. This proves  that
$I_{B_{R}}(Du^{n}- Du^{0})\to 0$ in 
$L_{\sfp,\sfq}(\bR^{d}_{S})$ and certainly in measure.
The theorem is proved. \qed

For integers $m\geq 1$ and $t_{0}\in(0,\infty)$  
introduce 
$$
\Gamma^{m}_{t_{0}}=\{(t_{1},...,t_{m}):t_{0}>
t_{1}>...>t_{m}>0=t_{m+1}\}.
$$
Also take Borel uniformly bounded functions
$\sigma^{ik}_{n}(s,x)$, $i=1,...,d$, $k=1,...,m$,
$n=0,1,2,...$, and assume that $\sigma^{ik}_{n}(s,x)
\to \sigma^{ik}_{0}(s,x)$ as $n\to\infty$ for almost all $(s,x)$.

Introduce  
\begin{equation}
                                                         \label{6.26.5}
Q^{k,n}_{s,t }f(x)=\sigma^{ik}_{n}(s,x)D_{i}T^{n}_{s,t }f(x).
\end{equation}

\begin{lemma}
                           \label{lemma 5.27.1}
Let $\phi_{t_{1},...,t_{m}}(x)$
be a bounded Borel function
on $\Gamma^{m}_{t_{0}}\times \bR^{d}$ such that
it equals zero whenever $\min(t_{k}-t_{k+1}:k=0,...,m )\leq\varepsilon$ for some $\varepsilon>0$, and
let $f\in L_{ p_{0}}$. Then 
$$  
\int_{\Gamma^{m}_{t_{0}}}T^{n}_{0,t_{m }}\big[
\phi_{t_{1},...,t_{m}}
Q^{m,n}_{t_{m},t_{m-1}}\cdot...\cdot
Q^{1,n}_{ t_{1},t_{0}}f \big](0)\,d t_{m} \cdot...\cdot d t_{1}
$$
\begin{equation}
                                 \label{5.27.6}
\to
\int_{\Gamma^{m}_{t_{0}}}T^{0} _{0,t_{m }}\big[
\phi_{t_{1},...,t_{m}}
Q^{m,0}_{t_{m},t_{m-1}}\cdot...\cdot
Q^{1,0}_{ t_{1},t_{0}}f \big](0)\,d t_{m} \cdot...\cdot d t_{1}
\end{equation}
as $n\to\infty$.
\end{lemma}

Proof. Set 
\begin{equation}
                                   \label{5.30.1}
\Gamma^{m}_{t_{0}}(\varepsilon)=
\Gamma^{m}_{t_{0}}\cap\{(t_{1},...,t_{m}):\min(t_{k}-t_{k+1}:k=0,...,m )\geq\varepsilon\}.
\end{equation}
Observe that on $\Gamma^{m}_{t_{0}}(\varepsilon)$ the functions
\begin{equation}
                                  \label{5.28.1}
T^{n}_{0,t_{m }}\big[
\phi_{t_{1},...,t_{m}}
Q^{m,n}_{t_{m},t_{m-1}}\cdot...\cdot
Q^{1,n}_{ t_{1},t_{0}}f \big](0)
\end{equation}
are uniformly bounded because the norms in $L_{ p_{0}}$
of the operators
 $Q^{m,n}_{t_{k},t_{k-1}}$ are uniformly bounded
and the $L_{ p_{0}}$-norms of the functionals
$T^{n}_{0,t_{m }}g(0) $ are uniformly bounded.

Next, we claim that the integrands in \eqref{5.27.6}
converge in measure on $\Gamma^{m}_{t_{0}}(\varepsilon)$. To prove the claim note that
with a constant $N$ independent of the $t_{k}$'s and $n$
$$
\big|T^{n}_{0,t_{m }}\big[
\phi_{t_{1},...,t_{m}}
Q^{m,n}_{t_{m},t_{m-1}}\cdot...\cdot
\{Q^{1,n}_{ t_{1},t_{0}}-Q^{1,0 }_{ t_{1},t_{0}}\}f \big](0)\big|
$$
$$
\leq N\|\{Q^{1,n}_{ t_{1},t_{0}}-Q^{1,0 }_{ t_{1},t_{0}}\}f\|_{L_{ p_{0}}}\to 0 
$$
in measure owing to Theorem \ref{theorem 5.22.1}.
It follows that while proving our claim we can
reduce the number of the operators $Q$ in \eqref{5.28.1} by one (replacing $f$ with $
Q^{1,0 }_{ t_{1},t_{0}} f$). By doing the same reduction
$m-1$ more times we see the it suffices to 
prove that $T^{n}_{0,t_{m }}\big[
\phi_{t_{1},...,t_{m}}f\big](0)\to
T^{0}_{0,t_{m }}\big[
\phi_{t_{1},...,t_{m}}f\big](0)$ in measure
on $\Gamma^{m}_{t_{0}}(\varepsilon)$.
But this holds even pointwise owing to Lemma
\ref{lemma 5.17.1}. This proves our claim
and with it the lemma as well. \qed

\mysection{A criterion for strong solutions
and the proof of Theorems \protect\ref{theorem 5.30.1} and  \protect\ref{theorem 6.1.1}}

                            \label{section 6.4.2}

We suppose that Assumptions \ref{assumption 12.19.1},
\ref{assumption 3.14.1}, and 
\ref{assumption 7.21.2}  are satisfied and  use the evolution family $T_{s,r}$ from
Section \ref{section 6.4.1}. We also use the notation (cf. 
\eqref{6.26.5})
$$
Q^{k }_{t,s}f(x)=\sigma^{ik} (t,x)D_{i}T _{t,s}f(x).
$$
Without loss of generality we concentrate on solutions
of \eqref{6.15.2} starting from the origin
and take $x_{t}$ as the trajectory of the Markov
process $X$ from Section \ref{section 6.4.1}
starting from the origin. As we know, it satisfies
\eqref{6.15.2} on a probability space 
$(\Omega,\cF,P)$ carrying
a $d_{1}$-dimensional Wiener process $w_{t}$. Recall that
for $x_{t}$ estimate \eqref{12.11.06} is available.

Just in case, keep in mind that $T_{s,r}f^{2}
=T_{s,r}(f^{2})$.

\begin{theorem}
                       \label{theorem 5.28.1}
 For each $f\in L_{ p_{0}}\cap L_{2 p_{0}}$ and $t_{0}>0$, for our solution $x_{s}$  of
  \eqref{6.15.2}, 
with  probability one,   we
have
\begin{equation}
                                                     \label{6.16.3}
 f(x_{t_{0}})=   T_{0,t_{0}}f(0)+\int_{0}^{t_{0}}
Q^{k }_{t_{1},t_{0}}f(x_{t_{1}}) \,dw^{k}_{t_{1}},
\end{equation}
\begin{equation}
                                                     \label{6.17.1}
   T_{0,t_{0}}f^{2}(0)=(   T_{0,t_{0}}f(0))^{2}+
\sum_{k}\int_{0}^{t_{0}}   T_{0,t_{1}}\big(
Q^{k }_{t_{1},t_{0}}f\big)^{2} (0)\,dt_{1}.
\end{equation}

\end{theorem}  

Proof. If $f\in C^{\infty}_{0}$, then $u(t,x):=  T_{t,t_{0}}f(x)$, $t\leq t_{0}$ is in $E^{1,2}_{\sfp,\sfq,\beta_{0}}(\bR^{d}_{t_{0}})$ by Lemma \ref{lemma 4.30.2}. Also $\cL u=0$ in $\bR^{d}_{t_{0}}$, which 
by It\^o's formula (Lemma 3.6  of \cite{13}) yields
that for any $R>0$
\begin{equation}
                                \label{5.28.4}
u(t_{0}\wedge \tau_{R},x_{t_{0}\wedge \tau_{R}})=u(0,0)+\int_{0}^{t_{0}\wedge \tau_{R}}
\sigma^{ik}(t_{1},x_{t_{1}})D_{i}u(t_{1},x_{t_{1}})
\,dw^{k}_{t_{1}},
\end{equation}
where $\tau_{R}$ is the first exit time of $(s,x_{s})$ from $C_{R}$. In addition, the stochastic integral
in \eqref{5.28.4} is a square integrable martingale, which owing to the boundedness of $u$ implies that
$$
\sum_{k}E\int_{0}^{t_{0}}\big[
\sigma^{ik}(t_{1},x_{t_{1}})D_{i}u(t_{1},x_{t_{1}})\big]^{2}\,dt_{1}<\infty.
$$
In turn, this allows us to let $R\to\infty$ in
\eqref{5.28.4} and obtain what is equivalent to
\eqref{6.16.3}. By taking the expectations of
the squares of both parts of \eqref{6.16.3}
we come to \eqref{6.17.1}.
This proves the theorem for $f\in C^{\infty}_{0}$.

In the general case observe that, as $f^{n}\in C^{\infty}_{0}$ tend to $f$ in $L_{ p_{0}}\cap L_{2 p_{0}}$,
we have $T_{0,t_{0}}f^{n}(0)\to T_{0,t_{0}}f(0)$ and
$T_{0,t_{0}}[f^{n}]^{2}(0)\to T_{0,t_{0}}f^{2}(0)$.
Furthermore, 
$\|Q^{k} _{t_{0},t_{1}}f^{n}-Q^{k}
   _{t_{0},t_{1}}f\|_{L_{ p_{0}}}\to 0$ for any  $t\in (0,t_{0})$. It follows that there is a subsequence $n'\to\infty$ such that
$$
\big(Q^{k}
   _{t_{0},t_{1}}f(x)\big)^{2}\leq\nliminf_{n'\to\infty}
\big(Q^{k}
   _{t_{0},t_{1}}f^{n'}(x)\big)^{2}
$$
for almost any $x\in\bR^{d}$ and   $t\in (0,t_{0})$.
Then Fatou's lemma and \eqref{6.17.1} allow us
to conclude that
\begin{equation}
                                                     \label{6.17.2}
   T_{0,t_{0}}f^{2}(0)\geq (   T_{0,t_{0}}f(0))^{2}+
\sum_{k}\int_{0}^{t}   T_{0,t_{1}}\big(
Q^{k }_{t_{1},t_{0}}f\big)^{2} (0)\,dt_{1}.
\end{equation}

Hence, the right-hand side of \eqref{6.16.3}
is well defined. Furthermore,
$$
E \Big|\int_{0}^{t_{0}}Q^{k }_{t_{1},t_{0}}f(x_{t_{1}}) \,dw^{k}_{t_{1}}-
\int_{0}^{t_{0}}\sigma^{ik}Q^{k }_{t_{1},t_{0}}f^{n}(x_{t_{1}}) \,dw^{k}_{t_{1}}\Big|^{2}
$$
$$
=\sum_{k}\int_{0}^{t_{0}}   T_{0,t_{1}} \big(Q^{k }_{t_{1},t_{0}}(f-f^{n})\big)^{2} (0)\,ds
$$
$$
\leq    T_{0,t_{0}}(f-f^{n})^{2}(0)- (   T_{0,t_{0}}(f-f^{n})(0))^{2},
$$
where the inequality is due to \eqref{6.17.2}. The last expression  
tends to zero, which allows us to get \eqref{6.16.3}
by passing to the limit in its version with $f^{n}$
in place of $f$. After that \eqref{6.17.1} follows as above.
The theorem is proved. \qed

By  Lemma \ref{lemma 4.28.1} and Theorem \ref{theorem 5.29.1} we have that
$Q^{k }_{t_{1},t_{0}}f\in L_{ p_{0}}\cap L_{2 p_{0}}$
for almost all $t_{1}<t_{0}$. Therefore,
by Theorem \ref{theorem 5.28.1} for almost all $t_{1}<t_{0}$
\begin{equation}
                                                     \label{6.17.20}
Q^{k_{1} }_{t_{1},t_{0}}f(x_{t_{1}})
=T_{0,t_{1}}Q^{k_{1} }_{t_{1},t_{0}}f(0)
+\int_{0}^{t_{1}}Q^{k_{2}}_{t_{2},t_{1}}
Q^{k_{1} }_{t_{1},t_{0}}f(x_{t_{2}})\,dw^{k_{2}}_{t_{2}}.
\end{equation}

After that we   substitute the result into \eqref{6.16.3}
to get
$$
f(x_{t_{0}})=T_{0,t_{0}}f(0) + \int_{0}^{t_{0}}
T_{0,t_{1}}Q^{k_{1} }_{t_{1},t_{0}}f(0) \,dw^{k_{1}}_{t_{1}}
$$  
\begin{equation} 
                              \label{6.17.3}
+\int_{0}^{t_{0}}\int_{0}^{t_{1}}
Q^{k_{2}}_{t_{2},t_{1}}
Q^{k_{1} }_{t_{1},t_{0}}f(x_{t_{2}})\,dw^{k_{2}}_{t_{2}} dw^{k_{1}}_{t_{1}}.
\end{equation}
Next, we apply Theorem \ref{theorem 5.28.1} to
$Q^{k_{2}}_{t_{2},t_{1}}
Q^{k_{1} }_{t_{1},t_{0}}f(x_{t_{2}})$ and get that
for almost all $t_{2}<t_{1}<t_{0}$
$$
Q^{k_{2}}_{t_{2},t_{1}}
Q^{k_{1} }_{t_{1},t_{0}}f(x_{t_{2}})=T_{0,t_{2}}
Q^{k_{2}}_{t_{2},t_{1}}
Q^{k_{1} }_{t_{1},t_{0}}f(0)+\int_{0}^{t_{2}}
Q^{k_{3}}_{t_{3},t_{2}}Q^{k_{2}}_{t_{2},t_{1}}
Q^{k_{1} }_{t_{1},t_{0}}f(x_{t_{3}})\,dw^{k_{3}}_{t_{3}}.
$$
After substituting the result in \eqref{6.17.3} we obtain
$$
f(x_{t_{0}})=T_{0,t_{0}}f(0) + \int_{0}^{t_{0}}
T_{0,t_{1}}Q^{k_{1} }_{t_{1},t_{0}}f(0) \,dw^{k_{1}}_{t_{1}}
$$  
\begin{equation} 
                              \label{6.10.1}
+\int_{0}^{t_{0}}\Big(\int_{0}^{t_{1}}
T_{0,t_{2}}
Q^{k_{2}}_{t_{2},t_{1}}
Q^{k_{1} }_{t_{1},t_{0}}f(0)\,dw^{k_{2}}_{t_{2}}\Big)\, dw^{k_{1}}_{t_{1}}
\end{equation}
$$
+\int_{0}^{t_{0}}\int_{0}^{t_{1}}\int_{0}^{t_{2}}
Q^{k_{3}}_{t_{3},t_{2}}Q^{k_{2}}_{t_{2},t_{1}}
Q^{k_{1} }_{t_{1},t_{0}}f(x_{t_{3}})\,dw^{k_{3}}_{t_{3}}dw^{k_{2}}_{t_{2}} dw^{k_{1}}_{t_{1}}.
$$

The formal objection to do these manipulations is that we should know
that the expression inside the parentheses in
\eqref{6.10.1} is, for instance,
predictable  as a function of $(\omega,t_{1})$ and this may not
happen if we allow any version of the inside stochastic integral
to be taken for each $t_{1}$.
A justification of the possibility to choose
right versions can be found either
in Remark 5.2 of \cite{Kr_21}
or in \cite{VK_76}.

By induction we obtain that for any  $t_{0}>0$ and $m\geq1$ (a.s.)  
$$
f(x_{t_{0}})=T_{0,t_{0}}f(0) 
$$
$$
+\sum_{n=1}^{m}\int_{\Gamma_{t_{0}}^{n}}T_{0,t_{n}}
Q^{k_{n}}_{t_{n},t_{n-1} }\cdot...\cdot
Q^{k_{1}}_{t_{1}, t_{0}}f(0)\,dw^{k_{n}}_{t_{n}}\cdot...\cdot dw^{k_{1}}_{t_{1}}
$$
\begin{equation}
                                                                \label{6.18.1}
+\int_{\Gamma_{t_{0}}^{m+1}}
Q^{k_{m+1}}_{t_{m+1},t_{m}}\cdot...\cdot
Q^{k_{1}}_{t_{1}, t_{0}}f(x_{t_{m+1}})\,dw^{k_{m+1}}_{t_{m+1}}\cdot...\cdot dw^{k_{1}}_{t_{1}},
\end{equation}
where by the expressions like 
$$
\int_{\Gamma_{t_{0}}^{m}}:::\,dw^{k_{m}}_{t_{m}}\cdot...\cdot dw^{k_{1}}_{t_{1}}
$$
we mean
$$
\int_{ 0}^{t_{0}}\,dw^{k_{1}}_{t_{1}}\int_{  0}^{t_{1}}\,dw^{k_{2}}_{t_{2}}
...\int_{ 0}^{t_{m-1}}:::\,dw^{k_{m}}_{t_{m}}.
$$

Introduce
$\bW^{m}_{t_{0}}$ as the closed linear
subspace of $L_{2}(\Omega,\cF,P)$ generated
by constants if $m=0$, by the set of
constants and
$$
\int_{0}^{t_{0}}f(t)\,dw^{k}_{t}
$$
if $m=1$
 or, if $m\geq 2$, by the set of constants and
$$
\int_{\Gamma^{n}_{t_{0}} }
f(t_{n},...,t_{1})
\,dw^{k_{n}}_{t_{n}}...\,dw^{k_{1}}_{t_{1}},
$$
where $k$, $(k_{1},...,k_{n})$, $n\leq m$ are arbitrary and $f(t)$ and $f(t_{n},...,t_{1})$ are arbitrary Borel bounded functions of their arguments. The projection operator in $L_{2}
(\Omega,\cF,P)$ on $\bW^{m}_{t_{0}}$ we denote
by $\Pi^{m}_{t_{0}}$.  

One knows (It\^o's theorem)
that $\xi\in L_{2}
(\Omega,\cF,P)$ is measurable with respect to
$\cF^{w}_{t_{0}}$ (the completion of $\sigma(w_{s}:s\leq t_{0})$) iff 
$$
E|\xi-\Pi^{m}_{t_{0}}\xi|^{2}\to 0
$$
as $m\to\infty$. Since the last term in \eqref{6.18.1}
is orthogonal to $\bW^{m}_{t_{0}}$,
for $\xi=f(x_{t_{0}})$ we 
have
$$
\xi-\Pi^{m}_{t_{0}}\xi
=\int_{\Gamma_{t_{0}}^{m+1}}
Q^{k_{m+1}}_{t_{m+1},t_{m}}\cdot...\cdot
Q^{k_{1}}_{t_{1}, t_{0}}f(x_{t_{m+1}})\,dw^{k_{m+1}}_{t_{m+1}}\cdot...\cdot dw^{k_{1}}_{t_{1}},
$$
$$
E|\xi-\Pi^{m}_{t_{0}}\xi|^{2} 
$$
\begin{equation}
                                 \label{5.29.5}
=\sum_{k_{1},...,k_{m+1}} \int_{\Gamma^{m+1}_{t_{0}} }T_{0,t_{m +1}}\big[
Q^{k_{m+1 }}_{t_{m+1 },t_{m }}\cdot...\cdot
Q^{k_{1}}_{t_{1},t_{0}}f\big]^{2}(0)
\,d t_{m +1} \cdot...\cdot d t_{1}
\end{equation}
 which yields the following.
 
\begin{theorem}
                                             \label{theorem 6.18.2}
Let $f\in L_{ p_{0}}\cap L_{2 p_{0}}$, $t_{0}>0$. Then
$f(x_{t_{0}})$ is $\cF^{w}_{t_{0}}$-measurable iff
\begin{equation}
                                                                \label{6.18.10}
 \sum_{k_{1},...,k_{m}} \int_{\Gamma^{m}_{t_{0}} }T_{0,t_{m }}\big[
Q^{k_{m }}_{t_{m },t_{m-1}}\cdot...\cdot
Q^{k_{1}}_{t_{1},t_{0}}f\big]^{2}(0)
\,d t_{m } \cdot...\cdot d t_{1}\to0 
\end{equation}
as $m\to\infty$.
Furthermore, under either of the above equivalent conditions
$$
f(x_{t_{0}})=T_{0,t_{0}}f(0) 
$$
\begin{equation}
                                                                \label{7.9.1}
+\sum_{n=1}^{\infty}\int_{\Gamma_{t_{0}}^{n}}T_{0,t_{n}}
Q^{k_{n}}_{t_{n},t_{n-1} }\cdot...\cdot
Q^{k_{1}}_{t_{1}, t_{0}}f(0)\,dw^{k_{n}}_{t_{n}}\cdot...\cdot dw^{k_{1}}_{t_{1}},
\end{equation}
where the series converges in the mean square sense.

\end{theorem}

Note that we obtained Theorem \ref{theorem 6.18.2}
only for special solutions $x_{t}$ and  not for any solution.

For the future observe that the right-hand side of
\eqref{5.29.5} is finite and is a decreasing function
of $m$.

 Next, we need a special approximations of $\sigma^{k}$
with functions that are smooth in $x$. Any smooth approximations of $b$ would suffice. Take   nonnegative $\xi\in C^{\infty}_{0}(\bR)$, $\eta\in C^{\infty}_{0}(\bR^{d})$
with unit integrals and supports in unit balls, set $\zeta(t,x)=\xi (t)\eta (x)$, $\zeta_{n}(t,x)=n^{d+1}\zeta(nt,nx)$,
$\eta_{n}( x)=n^{d }\eta( nx)$ and  define
 $$
\sigma^{(n)}=\sigma * \eta_{n}, 
\quad D\sigma^{(n)}_{M}=(D\sigma_{M}) * \eta_{n}, \quad
D\sigma^{(n)}_{B}=(D\sigma_{B}) * \eta_{n}, 
$$
where the convolution is performed with respect to $x$, and 
$$
b^{(n)}=b * \zeta_{n},\quad
b^{(n)}_{M}=b_{M} * \zeta_{n},\quad
b^{(n)}_{B}=b_{B} * \zeta_{n},
$$
where the convolution is performed with respect to $(t,x)$.
As is easy to see, for each $n$, $\sigma^{(n)}$
is bounded along with  any of its derivatives of any order with respect to $x$, $b^{(n)}$
 is a smooth bounded function.
Furthermore, Minkowski's inequality easily
shows that $D\sigma^{(n)},b^{(n)}$ (with the same  $p_{D\sigma},p_{b}$, $\rho_{D\sigma}$, $\rho_{b}$,
$r_{b}$, 
$\frp_{b},\frq_{b}$) satisfy
$$
\widehat{D\sigma^{(n)}_{M}}\leq \widehat{D\sigma_{M}},\quad  \overline{D\sigma^{n}}\leq \overline{D\sigma},\quad
\hat b^{(n)}_{M}
\leq \hat b_{M}, \quad \beta_{D\sigma^{(n)}}
\leq \beta_{D\sigma} .
$$

Finally, introduce  
$$
\Theta_{m}=\{t: \overline{D\sigma}(t)\leq
m\},
$$
$$
\sigma^{(n)}_{m}(t,x)=\sigma^{(n)} (t,x)
I_{\Theta_{m}}(t)+\kappa I_{\Theta^{c}_{m}}(t),
$$
where $\kappa$ is any fixed $d\times d_{1}$-matrix such that $\kappa \kappa^{*}=(\delta^{ij})$.

The following is a consequence of Lemma 5.11 of \cite{10}. 

\begin{lemma}
                    \label{lemma 6.13.1}
Set $a^{n}_{m}=\sigma^{(n)}_{m}\sigma^{(n)*}_{m}$. Then there is a sequence $m(n)\to\infty$
as $n\to\infty$ such that
for sufficiently large $n$ the eigenvalues of
$a^{n}_{m(n)}$ are between $\delta/4$ and $4/\delta$
if $N(d )\widehat{D\sigma}_{M}\leq 
\delta^{1/2}/4$, where $N(d)$ is specified
in the proof of Lemma 5.11 of \cite{10}.

\end{lemma}

It follows from the above that there exist sequences
  $b^{n},\sigma^{n}$ such that with $a^{n}=
\sigma^{n}(\sigma^{n}) ^{*}$ Assumption
\ref{assumption 5.16.1} is satisfied. In addition,
$\sigma^{n}\to \sigma=:\sigma^{0}$ (a.e.) as $n\to\infty$
and, for each $n>0$, the functions $b^{n},\sigma^{n}$ are smooth in $x$ with each derivative bounded with respect to $(t,x)$. Then we define $x^{n}_{s}$ as
(unique)  solutions of \eqref{6.15.2}
with $\sigma^{n}$ and $b^{n}$ in place of $\sigma$
and $b$, respectively, and with $(0,0)$ as the initial data. These processes satisfy \eqref{5.28.2}
with $\cSigma^{n}=(\sigma^{n}(\sigma^{n})^{*})^{1/2}$,
$b^{n}$ in place of $\cSigma$, $b$
and, as it is pointed out after Assumption 
\ref{assumption 5.16.1}, converge in distribution
to a process $x_{t}$, which satisfies (on a probability space) equation
\eqref{5.28.2} with $\cSigma =(\sigma  \sigma  ^{*})^{1/2}$ and zero initial data and, therefore, satisfies
\eqref{6.15.2} on an extended probability space.
Since $x^{n}_{s}$ are strong solutions, we may assume
that all of them are given on the same probability space
as $x_{s}$ 
with the same $d_{1}$-dimensional Wiener process.

Then we can apply the results of Section \ref{section 6.4.1}
given after Assumption
\ref{assumption 5.16.1} to the distributions of 
$x_{t}$   and $x^{n}_{s}$. By Theorem  5.10
of \cite{10} (here we use that $ p_{0}\in ( d/2  , p_{b})$),
for $f\in L_{ p_{0}}\cap L_{2 p_{0}}$ and $\xi_{n}
=f(x^{n}_{t_ 0})$ we have
\begin{equation}
                            \label{6.11.4}
\sup_{n>0}\sum_{m=1}^{\infty}E|\xi_{n}-\Pi^{m}_{t_{0}}\xi_{n}|^{2}
\leq N\|f\|_{L_{2p_{0}}},
\end{equation}
where  $N$ depends only on $d,\delta,p_{D\sigma},p_{b},p_{0},
\beta_{D\sigma}(\infty),\beta_{b}(\infty),\rho_{\sigma},\rho_{b},t_{0}$
 (independent of $r_{b}$).
 
{\bf Proof of Theorems \ref{theorem 5.30.1}
and \ref{theorem 6.1.1}}. 
To prove Theorem \ref{theorem 5.30.1}, we take $x_{t}$ from above and observe that to show that it is a strong solution it suffices to show
that the series composed of the terms in \eqref{6.18.10} converges whenever $t_{0}>0$
and $f\in L_{ p_{0}}\cap L_{2 p_{0}}$. In light of
\eqref{6.11.4} and Fatou's lemma,
to prove that,
it suffices to show that for each $m$
$$
E|f(x_{t_{0}}) -\Pi^{m-1}_{t_{0}}f(x_{t_{0}})|^{2}
\leq\nliminf_{n\to\infty}E|f(x^{n}_{t_{0}})-\Pi^{m-1}_{t_{0}}f(x^{n}_{t_{0}})|^{2}
$$
or, which is slightly stronger (see \eqref{5.29.5}), that
for each $k_{1},...,k_{m}$ 
$$
\int_{\Gamma^{m}_{t_{0}} }T_{0,t_{m }}\big[
Q^{k_{m }}_{t_{m },t_{m-1}}\cdot...\cdot
Q^{k_{1}}_{t_{1},t_{0}}f\big]^{2}(0)
\,d t_{m } \cdot...\cdot d t_{1}
$$
\begin{equation}
                                  \label{5.30.2}
\leq \nliminf_{n\to\infty}
\int_{\Gamma^{m}_{t_{0}} }T^{n}_{0,t_{m }}\big[
Q^{k_{m }n}_{t_{m },t_{m-1}}\cdot...\cdot
Q^{k_{1}n}_{t_{1},t_{0}}f\big]^{2}(0)
\,d t_{m } \cdot...\cdot d t_{1}=:I.
\end{equation}

Fix $t_{0}>0$, $m$, and $k_{1},...,k_{m}$. Take  $\Gamma^{m}_{t_{0}}(\varepsilon)
$ from \eqref{5.30.1} and for $\lambda >0$ set
$$
\Lambda (\varepsilon,\lambda)
=\{(t_{1},...,t_{m},x):(t_{1},...,t_{m})
\in \Gamma^{m}_{t_{0}}(\varepsilon),
|Q^{k_{m }}_{t_{m },t_{m-1}}\cdot...\cdot
Q^{k_{1}}_{t_{1},t_{0}}f(x)|\leq \lambda\},
$$
$$
\phi^{\varepsilon,\lambda}_{t_{1},...,t_{m}}(x)=I_{\Lambda (\varepsilon,\lambda)}
Q^{k_{m }}_{t_{m },t_{m-1}}\cdot...\cdot
Q^{k_{1}}_{t_{1},t_{0}}f(x) .
$$

By Lemma  \ref{lemma 5.27.1} and H\"older's
inequality
$$
\int_{\Gamma^{m}_{t_{0}}}T^{0} _{0,t_{m }}\big[
\phi^{\varepsilon,\lambda}_{t_{1},...,t_{m}}
Q^{k_{m} }_{t_{m},t_{m-1}}\cdot...\cdot
Q^{k_{1}}_{ t_{1},t_{0}}f \big](0)\,d t_{m} \cdot...\cdot d t_{1}
$$
$$ 
=\lim_{n\to\infty}\int_{\Gamma^{m}_{t_{0}}}T^{n}_{0,t_{m }}\big[
\phi^{\varepsilon,\lambda}_{t_{1},...,t_{m}}
Q^{k_{m},n}_{t_{m},t_{m-1}}\cdot...\cdot
Q^{k_{1},n}_{ t_{1},t_{0}}f \big](0)\,d t_{m} \cdot...\cdot d t_{1}
$$
\begin{equation}
                                 \label{5.27.06}
\leq I^{1/2}\Big(\lim _{n\to\infty}
\int_{\Gamma^{m}_{t_{0}}}T^{n} _{0,t_{m }}\big(\big[
\phi^{\varepsilon,\lambda}_{t_{1},...,t_{m}}
 \big]^{2}\big)(0)\,d t_{m} \cdot...\cdot d t_{1}
\Big)^{1/2}.
\end{equation}

Here the function $\big[
\phi^{\varepsilon,\lambda}_{t_{1},...,t_{m}}
 \big]^{2}$ is bounded and
by Lemma \ref{lemma 5.17.1}
$$
T^{n} _{0,t_{m }}\big(\big[
\phi^{\varepsilon,\lambda}_{t_{1},...,t_{m}}
 \big]^{2}\big)(0)\to
T  _{0,t_{m }}\big(\big[
\phi^{\varepsilon,\lambda}_{t_{1},...,t_{m}}
 \big]^{2}\big)(0)
$$
for any $t_{1},...,t_{m}$. It follows that
$$
\int_{\Gamma^{m}_{t_{0}} }T_{0,t_{m }}I_{\Lambda (\varepsilon,\lambda)}\big[
Q^{k_{m }}_{t_{m },t_{m-1}}\cdot...\cdot
Q^{k_{1}}_{t_{1},t_{0}}f\big]^{2}(0)
\,d t_{m } \cdot...\cdot d t_{1}
$$
$$
=\int_{\Gamma^{m}_{t_{0}}}T  _{0,t_{m }}\big[
\phi^{\varepsilon,\lambda}_{t_{1},...,t_{m}}
Q^{k_{m}}_{t_{m},t_{m-1}}\cdot...\cdot
Q^{k_{1} }_{ t_{1},t_{0}}f \big](0)\,d t_{m} \cdot...\cdot d t_{1}
$$
$$
\leq I^{1/2}\Big(
\int_{\Gamma^{m}_{t_{0}} }T_{0,t_{m }}\big[
Q^{k_{m }}_{t_{m },t_{m-1}}\cdot...\cdot
Q^{k_{1}}_{t_{1},t_{0}}f\big]^{2}(0)
\,d t_{m } \cdot...\cdot d t_{1}\Big)^{1/2}.
$$
By sending $\varepsilon\downarrow0$ and $
\lambda \to\infty$ in the first expression
we come to what is equivalent to \eqref{5.30.2}.

Thus indeed
$x_{t}$ is a strong solution.

That  estimate  \eqref{5.31.1}   holds follows from what is said after
Assumption \ref{assumption 5.16.1} and
Lemma \ref{lemma 6.1.2}, which also proves
Theorem \ref{theorem 6.1.1}.   

 To prove the second assertion of Theorem \ref{theorem 5.30.1} observe that  by Theorem 4.4 of \cite{10} all solutions
of \eqref{6.15.2} possessing property ($\text{a}\,'$)  
and having the same initial data have the same
distributions. Then the second assertion  
follows directly from  Theorem 5.1 of \cite{10}.
\qed

From the above proof we have the following.

\begin{corollary}
                        \label{corollary 7.26.1}
For any $f\in L_{p_{0}}\cap L_{2p_{0}}$, $t_{0}>0$,
and $\xi:=f(x_{t_{0}})$ we have
$$
\sum_{m=1}^{\infty}E|\xi -\Pi^{m}_{t_{0}}\xi |^{2}
\leq N\|f\|_{L_{2p_{0}}},
$$
where $N$ is taken from \eqref{6.11.4}.
\end{corollary}

\mysection{Proof of Theorem \protect\ref{theorem 7.18.1}}  
                           \label{section 7.29.1}

Observe that  by Theorem 4.4 of \cite{10} all solutions
of \eqref{6.15.2}  satisfying \eqref{7.19.1} and 
 having the same initial data have the same
distributions. Then  the uniqueness statement in 
Theorem \ref{theorem 7.18.1} 
follows directly from  Theorem 5.1 of \cite{10},
provided that there exists a strong solution
possessing property (a). Therefore, we now concentrate
on the existence of such solutions.

First we prove an auxiliary fact which is probably well known but it is easier to prove it than to find it in the literature.

\begin{lemma}
                          \label{lemma 7.27.1}
For any $m,n=1,2,...$, $t>0$, and 
$$
f(t_{1},...,t_{m})=\{
f^{k_{m},...,k_{1}}(
t_{1},...,t_{m}),k_{i}=1,...,d_{1}\},
$$
 given
on $\Gamma^{m}_{t}$ and square integrable there
we have
\begin{equation}
                         \label{7.27.1}
E\Big(\int_{\Gamma^{m}_{t}}f(
t_{1},...,t_{m})\,dw_{t_{m}}...dw_{t_{1}}\Big)^{2n}
\leq N \|f\|_{L_{2}(\Gamma^{m}_{t})}^{2n},
\end{equation}
where $N\, (<\infty)$ depends only $n,m,d_{1}$, and by the repeated stochastic
integral above we mean
\begin{equation}
                         \label{7.27.3}
\sum_{k_{1},...,k_{m}}
\int_{\Gamma^{m}_{t}}f^{k_{m},...,k_{1}}(
t_{1},...,t_{m})\,dw^{k_{m}}_{t_{m}}...dw^{k_{1}}_{t_{1}}
\end{equation}
\end{lemma}

Proof. Clearly, it suffices to prove
\eqref{7.27.1} for each particular term in
\eqref{7.27.3}. Introduce $A_{n,m}$
as the supremum of
$$
E\Big(\int_{\Gamma^{m}_{t}}f(
t_{1},...,t_{m})\,dB^{m}_{t_{m}}...dB^{1}_{t_{1}}\Big)^{2n}
$$
taken over all sets of
$\{B^{1}_{\cdot},...,
B^{m}_{\cdot}\}\subset \{w^{1}_{\cdot},...,w^{d_{1}}_{\cdot}\}$ and functions $f(
t_{1},...,t_{m})$ on $\Gamma^{m}_{t}$ having the $L_{2}$-norm equal to one. To prove the lemma, we 
only need
to show that $A_{n,m}<\infty$ for all $n,m$.

We are going to use the induction on $m$.
If $m=1$, the stochastic integral is normally distributed and \eqref{7.27.1} is obvious,
so $A_{n,1}<\infty$

Suppose that for some $m\geq1$ and any $n=1,2,...$
we have $A_{n,m}<\infty$.
Then take $f (t_{1},...,t_{m+1})$ such that
$\|f\|_{L_{2}(\Gamma^{m+1}_{t})}=1$, and observe
that by Burknolder-Davis-Gundy inequality
$$
I:=E\Big(\int_{\Gamma^{m+1}_{t}}f(
t_{1},...,t_{m})\,dB^{m+1}_{t_{m+1}}...dB^{1}_{t_{1}}\Big)^{2 n}\leq N(n)E\Big( 
\int_{0}^{t}I^{2}( t_{1})\,dt_{t}\Big)^{n}
$$
$$
=
N(n)\int_{(0,t)^{n}}EI^{2}( s_{1})\cdot...\cdot
I^{2}( s_{n})\,ds_{1}...ds_{n},
$$
where
$$
I (s)= 
\int_{\Gamma^{m}_{s}}f 
(s,t_{2},....,t_{m+1})\,dB^{m+1 }_{t_{m+1}}...
dB^{2}_{t_{2}}.
$$
Introduce $J(s)$ by $I(s)=J(s)\|f(s,\cdot)\|_{L_{2}(\Gamma^{m}_{s})}$ and observe that
$$
EJ^{2}(s_{1})\cdot...\cdot J^{2}(s_{n})
\leq \Big(\prod_{i=1}^{n} EJ^{2n}(s_{i})\Big)^{1/n}
\leq A_{n,m},
$$
where the last inequality holds by assumption. It follows that
$$
I\leq N(n)A_{n,m}
\int_{(0,t)^{n}}\|f(s_{1},\cdot)\|^{2}_{L_{2}(\Gamma^{m}_{s})}\cdot...\cdot
\|f(s_{1},\cdot)\|^{2}_{L_{2}(\Gamma^{m}_{s_{n}})}\,ds_{1}...ds_{n}.
$$
Since the last integral is, obviously,
$\|f\|_{L_{2}(\Gamma^{m+1}_{t})}^{2n}=1$, we have
$I\leq N(n) A_{n,m}$ and the arbitrariness
of $f$ and $B^{i}_{\cdot}$ implies that $A_{n,m+1}\leq
N(n) A_{n,m}$. This proves the lemma. \qed
 
Next, we specify
  $\varepsilon_{0},\varepsilon_{0}'$ 
in Theorem \ref{theorem 7.18.1} 
in the following way. 
Define $\frp_{b}=p_{b},\frq_{b}=q_{0}+1$.    Then set 
$\varepsilon_{0}
=\varepsilon_{\sigma}$, where $\varepsilon_{\sigma}$
is taken from Assumption \ref{assumption 7.21.2}
and reduce $\rho_{\sigma}$, if necessary
(depending on the modulus of continuity of $\beta_{D\sigma}$) so that 
$\beta_{D\sigma}(\rho^{2}_{\sigma})\leq \varepsilon_{\sigma}$.
Finally,  set $  \varepsilon'_{0}=\varepsilon_{b}/2$, where $\varepsilon_{b}$
is taken from Assumption \ref{assumption 7.21.2}.
After that we suppose that 
$\widehat{{D\sigma_{M}}}\leq \varepsilon_{0},
 \hat b_{M}\leq \varepsilon'_{0}$, take
a solution $x_{t}$ of \eqref{6.15.2} with $(t,x)=(0,0)$
  possessing property (a),
which exists by Theorem 3.9 of \cite{10} and set up the  goal of proving that $x_{t}$ is a strong solution.

Let $(\Omega,\cF,P)$ be the probability space on which
$x_{t}$ is defined and let $w_{t}$ be the corresponding Wiener process. For $n=1,2,...$, define
$$
b_{n}=b_{M}+b_{B}I_{|b_{B}|\leq n},\quad
b_{nM}=b_{M},\quad b_{nB}=b_{B}I_{|b_{B}|\leq n},
$$
$$
\gamma_{n}=\sigma^{*}a^{-1}b_{B}I_{|b_{B}|> n},
\quad \phi_{n}=-\int_{0}^{\infty}\gamma_{n}(s,x_{s})\,dw_{s}-(1/2)\int_{0}^{\infty}|\gamma_{n}(s,x_{s})
|^{2}\,ds.
$$
Observe that $\gamma_{n}(t,x)$ is bounded by
a function of $t$ which is square integrable
over $(0,\infty)$ (see \eqref{7.26.1}). It follows
that for any $\alpha\in\bR$ we have
$E\exp(\alpha\phi_{n})<\infty$. Another useful fact
following from \eqref{7.26.1}
is that
\begin{equation}
                        \label{7.28.3}
\int_{\bR}\sup_{\bR^{d}}|\gamma_{n}(t,x)|^{2}\,dt \to 0
\end{equation}
as $n\to\infty$.

Next, introduce $P^{n}(d\omega)=e^{\phi_{n}}P(d\omega)$,
$$
w^{(n)}_{t}=w_{t}+\int_{0}^{t}\gamma_{n}(s,x_{s})\,ds.
$$
By Girsanov's theorem $P^{n}$ is a probability measure, $w^{(n)}_{t}$ is a Wiener process
on $(\Omega,\cF,P^{n})$, and
\begin{equation}
                                \label{7.26.2}
x_{t}=\int_{0}^{t}\sigma(s,x_{s})\,dw^{(n)}_{s}
+\int_{0}^{t}b_{n}(s,x_{s})\,ds.
\end{equation}

Observe that for any $r\leq (2n)^{-1}=:r^{n}_{b}$ and $C\in \bC_{r}$
$$
\dashnorm b_{n}\|_{L_{\frp_{b},\frq_{b}}(C)}
\leq \dashnorm b_{M}\|_{L_{p_{b},\frq_{b}}(C)}
+n \leq r^{-1}\varepsilon_{b}/2+n\leq r^{-1}\varepsilon_{b}.
$$
Furthermore,  by using H\"older's inequality and
property (a), one easily sees that for any $T\in(0,\infty)$,
$m=1,2,...$, and Borel nonnegative $f$
($E^{n}$ is the expectation sign corresponding to $P^{n}$)
$$
E^{n}\Big(\int_{0}^{T}f(s,x_{s})\,ds\big)^{m}
\leq N\|f\|_{L_{p_{0},q_{0}}},
$$
where $N$ is independent of $f$. Now the assertion
of uniqueness in Theorem \ref{theorem 5.30.1}
and Corollary \ref{corollary 7.26.1} imply
that for $f\in L_{p_{0}}\cap L_{2p_{0}}$, $t_{0}>0$ on $\Gamma^{m}_{t_{0}}$, $m=1,2,...$,
for $k_{i}=1,...,d_{1}$, $i=1,...,m$, there exist deterministic functions 
$f^{n,k_{m},...,k_{1}} (t_{1},...,t_{m})$ 
square integrable over $\Gamma^{m}_{t_{0}}$  such that
$$
\sum_{m=1}^{\infty}
E^{n}\Big|f(x_{t_{0}})-c_{n}
$$
\begin{equation}
                             \label{7.26.3}
-
\sum_{i=1}^{m}\sum_{k_{1},...,k_{i}}
\int_{\Gamma^{i}_{t_{0}}}
f^{n,k_{i},...,k_{1}} (t_{1},...,t_{i})
\,dw^{(n)k_{i}}_{t_{i}}\cdot...\cdot dw^{(n)k_{1}}_{t_{1}}
\Big|^{2}\leq N\|f\|_{L_{2p_{0}}},
\end{equation}
where $c_{n}=E^{n}f(x_{t_{0}})$ and $N$ is the constant from
\eqref{6.11.4}.

Since
$$
c_{n}^{2}+\sum_{i=1}^{\infty}\sum_{k_{1},...,k_{i}}
\|f^{n,k_{i},...,k_{1}}\|^{2}_{L_{2}(\Gamma^{i}_{t_{0}})}=E^{n}f^{2}(x_{t_{0}}),
$$
and the right-hand side is bounded by a constant
independent of $n$, there is a subsequence $n'\to\infty$ such that $f^{n',k_{i},...,k_{1}}$
converge weakly in $L_{2}(\Gamma^{i}_{t_{0}})$
to certain functions $f^{ k_{i},...,k_{1}}$.
Of course, $c_{n}\to Ef(x_{t_{0}})$.

Next we prove three auxiliary facts.

\begin{lemma}
                        \label{lemma 7.28.1}
If $f\in L_{2}(\Gamma^{m}_{t_{0}})$, then
for any $n,k=1,2,...$ and $k_{1},...,k_{m}\in\{1,...,d_{1}\}$
$$
I:=E\Big(\int_{\Gamma^{m}_{t_{0}}}
f(t_{1},...,t_{m})\,dw^{(n)k_{m}}_{t_{m}}\cdot...
\cdot dw^{(n)k_{1}}_{t_{1}}\Big)^{2k}\leq 
N\|f\|_{L_{2}(\Gamma^{m}_{t_{0}})}^{2k},
$$
where (note $E$ not $E^{n}$) $N$ depends only on $m,k,d,d_{1},\delta$,
and $\beta_{b}(\infty)$.

\end{lemma}

The proof of the lemma is obtained by observing that
owing to Lemma \ref{lemma 7.27.1} and Girsanov's theorem
$$
I=Ee^{\phi_{n}}\Big(\int_{\Gamma^{m}_{t_{0}}}
f(t_{1},...,t_{m})\,dw^{ k_{m}}_{t_{m}}\cdot...
\cdot dw^{ k_{1}}_{t_{1}}\Big)^{2k}
$$
$$
\leq \Big(Ee^{2\phi_{n}}\Big)^{1/2}
\Big(E\Big(\int_{\Gamma^{m}_{t_{0}}}
f(t_{1},...,t_{m})\,dw^{ k_{m}}_{t_{m}}\cdot...
\cdot dw^{ k_{1}}_{t_{1}}\Big)^{4k}\Big)^{1/2}.
$$ \qed

\begin{lemma}
                        \label{lemma 7.28.2}
If $f\in L_{2}(\Gamma^{m}_{t_{0}})$, then
for any $n=1,2,...$ and $k_{1},...,k_{m}\in\{1,...,d_{1}\}$
$$
E\Big(\int_{\Gamma^{m}_{t_{0}}}
f  (t_{1},...,t_{m})
\,dw^{(n) k_{m}}_{t_{m}}\cdot...\cdot dw^{ (n)k_{1}}_{t_{1}}-
\int_{\Gamma^{m}_{t_{0}}}
f  (t_{1},...,t_{m})
\,dw^{ k_{m}}_{t_{m}}\cdot...\cdot dw^{ k_{1}}_{t_{1}}\Big)^{2}
$$
\begin{equation}
                                 \label{7.28.4}
\leq \varepsilon_{n}\|f\|_{L_{2}(\Gamma^{m}_{t_{0}})}^{2},
\end{equation}
where $\varepsilon_{n}$ is independent of $f$
and $\varepsilon_{n}\to0$ as $n\to\infty$.
\end{lemma}

Proof. Having in mind a usual telescoping procedure
we see that it suffices to prove that for $i=0,...,m$,
with obvious agreements in the extreme cases $i=0$ or $m$,
(keep in mind that for $i\leq m$ we set $\prod_{j=m}^{i}dw_{t_{j}}=dw_{t_{m}}...dw_{t_{i}}$)
$$
K_{i,n}:=E\Big(\int_{\Gamma^{m}_{t_{0}}}
f (t_{1},...,t_{m})
\prod_{j=m}^{i}  
dw^{(n)k_{ j}}_{t_{ j}}dw_{t_{i-1}}^{k_{i-1}}
\cdot...\cdot dw^{k_{1}}_{t_{1}}
$$
\begin{equation}
                                         \label{7.28.04}
-\int_{\Gamma^{m}_{t_{0}}}
f (t_{1},...,t_{m})
\prod_{j=m}^{i+1}  
dw^{(n)k_{ j}}_{t_{ j}}dw_{t_{i }}^{k_{i }}
\cdot...\cdot dw^{k_{1}}_{t_{1}}\Big)^{2}\leq \varepsilon_{n}\|f\|_{L_{2}(\Gamma^{m}_{t_{0}})}^{2}
\end{equation}

 Note that
with $\gamma_{n}(  t_{i} ):=\gamma_{n}( t_{i} ,x_{t_{i}})$ we have
$$
K_{i,n}
=E\Big(\int_{\Gamma^{m}_{t_{0}}}
f (t_{1},...,t_{m})
\prod_{j=m}^{i+1}  
dw^{(n)k_{ j}}_{t_{ j}}  \gamma_{n}^{k_{i}} (  t_{i} )dt_{ i }dw_{t_{ i-1}}^{k_{i-1}}
\cdot...\cdot dw^{k_{1}}_{t_{1}}\Big)^{2} 
$$
 $$
=\int_{\Gamma^{i-1}_{t_{0}}}
J _{n}(t_{i-1},...,t_{1})
dt_{i-1}\cdot...\cdot dt_{1},
$$
where
$$
J _{n}(t_{i-1},...,t_{1})=E\Big(\int_{\Gamma_{t_{i-1        }}^{m-i+1}}
f (t_{1},...,t_{m})
\prod_{j=m}^{i+1}  
dw^{(n)k_{ j}}_{t_{ j}} \gamma_{n}^{k_{i}} (  t_{i} )dt_{ i }\Big)^{2}
$$
$$
\leq E\int_{0}^{t_{i-1}}|\gamma_{n}^{k_{i}}|^{2} (  t_{i} )dt_{ i }\int_{0}^{t_{i-1}}
\Big(\int_{\Gamma_{t_{i }}^{m-i }}
f (t_{1},...,t_{m})
\prod_{j=m}^{i+1}  
dw^{(n)k_{ j}}_{t_{ j}}\Big)^{2}dt_{i}.
$$
Here the first integral under the expectation sign
tends to zero as $n\to\infty$ uniformly
with respect to $t_{i-1},\omega$ (see \eqref{7.26.1}) and 
$$
E \int_{0}^{t_{i-1}}
\Big(\int_{\Gamma_{t_{i }}^{m-i }}
f (t_{1},...,t_{m})
\prod_{j=m}^{i+1}  
dw^{(n)k_{ j}}_{t_{ j}}\Big)^{2}dt_{i}
$$
$$
=\int_{\Gamma_{t_{i-1 }}^{m-i+1 }}| 
f (t_{1},...,t_{m})|^{2}dt_{m}\cdot...\cdot
dt_{i}.
$$
This easily implies \eqref{7.28.04} and the lemma is proved. \qed

\begin{lemma}
                                \label{lemma 7.29.1}
Let $f^{n}\to f$ weakly in $L_{2}(\Gamma^{m}_{t_{0}})$ 
as $n\to \infty$ and $k_{1},...,k_{m}\in\{1,...,d_{1}\}$. Then
$$
 \int_{\Gamma^{m}_{t_{0}}}
f^{n}  (t_{1},...,t_{m})
\,dw^{ k_{m}}_{t_{m}}\cdot...\cdot dw^{ k_{1}}_{t_{1}}
\to \int_{\Gamma^{m}_{t_{0}}}
f  (t_{1},...,t_{m})
\,dw^{ k_{m}}_{t_{m}}\cdot...\cdot dw^{ k_{1}}_{t_{1}}
$$
weakly in $L_{2}(\Omega)$ as $n\to \infty$.
\end{lemma}

To prove the lemma, it suffices to observe that
for any $\eta\in L_{2}(\Omega)$ the
functional
$$
E\eta\int_{\Gamma^{m}_{t_{0}}}
f  (t_{1},...,t_{m})
\,dw^{ k_{m}}_{t_{m}}\cdot...\cdot dw^{ k_{1}}_{t_{1}}
$$
is bounded in $ L_{2}(\Gamma^{m}_{t_{0}})$, hence
continuous and weakly continuous. \qed

Now note that

$$
M^{n'}:=\int_{\Gamma^{m}_{t_{0}}}
f^{n',k_{m},...,k_{1}} (t_{1},...,t_{m})
\,dw^{(n')k_{m}}_{t_{m}}\cdot...\cdot dw^{(n')k_{1}}_{t_{1}} 
$$
$$
=:\int_{\Gamma^{m}_{t_{0}}}
f^{n',k_{m},...,k_{1}} (t_{1},...,t_{m})
\,dw^{ k_{m}}_{t_{m}}\cdot...\cdot dw^{ k_{1}}_{t_{1}}
+J^{n'}=:I^{n'}+J^{n'},
$$
where $J^{n'}\to0$ in $L_{2}(\Omega)$ as $n'\to
\infty$ by Lemma \ref{lemma 7.28.2} and by Lemma \ref{lemma 7.29.1}
$$
I^{n'}\to
\int_{\Gamma^{m}_{t_{0}}}
f^{ k_{m},...,k_{1}} (t_{1},...,t_{m})
\,dw^{ k_{m}}_{t_{m}}\cdot...\cdot dw^{ k_{1}}_{t_{1}}
=:M
$$
weakly in $L_{2}(\Omega)$. Since $e^{\phi_{n}/2}\to1$ strongly
in $L_{2}(\Omega)$, we also have
that $e^{\phi_{n'}/2}M^{n'}\to M$
weakly in $L_{2}(\Omega)$.
Also, obviously, $e^{\phi_{n'}/2}
  f(x_{t_{0}})\to f(x_{t_{0}})$ and
$e^{\phi_{n'}/2} c_{n'}\to Ef(x_{t_{0}})$ weakly (strongly) in  $L_{2}(\Omega)$.

By Fatou's lemma the sum of the $\nliminf$'s of the
terms on the left-hand side of \eqref{7.26.3} with $n'$ in place of $n$
is less than the $\nliminf$ of the left-hand side
of \eqref{7.26.3} and, hence, is finite. Taking into account
that ``the norm of the weak limit is less than
the $\nliminf$ of the norms'' and taking into account
the above results we conclude that
$$
\sum_{m=1}^{\infty}
E^{n}\Big|f(x_{t_{0}})-Ef(x_{t_{0}})
$$
$$
-
\sum_{i=1}^{m}\sum_{k_{1},...,k_{i}}
\int_{\Gamma^{i}_{t_{0}}}
f^{ k_{i},...,k_{1}} (t_{1},...,t_{i})
\,dw^{(n)k_{i}}_{t_{i}}\cdot...\cdot dw^{(n)k_{1}}_{t_{1}}
\Big|^{2}\ <\infty,
$$
which implies that $f(x_{t_{0}})$ is $\cF^{w}_{t_{0}}$-measurable and the arbitrariness of $f$ and $t_{0}$,
finally, 
bring the proof of the theorem to an end.\qed

\end{document}